\documentclass[a4paper]{amsart}
\usepackage{amssymb}
\usepackage[shortlabels]{enumitem}
\usepackage{hyperref}
\usepackage[capitalize]{cleveref}
\usepackage[labelfont={}]{subcaption}
\usepackage{thm-restate}
\usepackage{microtype}

\usepackage{tikz}
\usepackage{quiver}
\usetikzlibrary{backgrounds}

\usepackage{float}

\newtheorem{theorem}{Theorem}[section]
\newtheorem{corollary}[theorem]{Corollary}
\newtheorem{lemma}[theorem]{Lemma}
\newtheorem{proposition}[theorem]{Proposition}
\newtheorem{question}[theorem]{Question}

\newtheorem{letterthm}{Theorem}

\crefname{letterthm}{Theorem}{Theorems}
\Crefname{letterthm}{Theorem}{Theorems}
\newtheorem{lettercor}[letterthm]{Corollary}

\crefname{lettercor}{Corollary}{Corollaries}
\Crefname{lettercor}{Corollary}{Corollaries}

\crefname{proposition}{Pro\-po\-si\-tion}{Pro\-po\-si\-tions}
\Crefname{proposition}{Pro\-po\-si\-tion}{Pro\-po\-si\-tions}

\theoremstyle{definition}
\newtheorem{definition}[theorem]{Definition}
\newtheorem{remark}[theorem]{Remark}
\newtheorem{example}[theorem]{Example}
\theoremstyle{remark}
\numberwithin{equation}{theorem}

\newcommand{\Dfn}[1]{\emph{\color{blue}#1}} 

\renewcommand{\mod}{\operatorname{mod}}
\newcommand{\findim}{\operatorname{findim}}
\newcommand{\proj}{\operatorname{proj}}
\newcommand{\Ext}{\operatorname{Ext}}
\newcommand{\End}{\operatorname{End}}
\newcommand{\Hom}{\operatorname{Hom}}
\newcommand{\sAus}{\operatorname{\underline{Aus}}}
\newcommand{\pstab}[1]{\underline{P}_{#1}}
\newcommand{\Z}{\mathbb{Z}}
\newcommand{\domdim}{\operatorname{domdim}}
\newcommand{\codomdim}{\operatorname{codomdim}}
\newcommand{\gordim}{\operatorname{Gordim}}
\newcommand{\fingordim}{\operatorname{finGordim}}
\newcommand{\height}{\operatorname{ht}}
\newcommand{\ceilfindim}[1]{\left\lceil\frac{\findim #1}{2}\right\rceil}
\newcommand{\idim}{\operatorname{idim}}
\newcommand{\pdim}{\operatorname{pdim}}
\newcommand{\gldim}{\operatorname{gldim}}
\newcommand{\DynA}{\mathbb{A}}
\NewDocumentCommand{\indec}{O{i}O{k}}{e_{#1}A / e_{#1}J^{#2}}
\newcommand{\dist}{\operatorname{dist}}
\newcommand{\level}{\operatorname{lvl}}
\NewDocumentCommand{\mik}{O{i}O{k}}{M(#1,#2)}
\newcommand{\kupisch}{\left[c_0, c_1,\dots, c_{n-1}\right]}
\newcommand{\cokupisch}{\left[d_0, d_1,\dots, d_{n-1}\right]}
\NewDocumentCommand{\cosyz}{O{}}{\ensuremath{\Sigma^{#1}}}
\newcommand{\psiforest}{F_\psi}
\newcommand{\gammaforest}{F_\gamma}
\newcommand{\psiquiver}{Q_\psi}
\NewDocumentCommand{\stabidim}{O{M}}{\idim\pstab{#1}}
\newcommand{\op}{\textnormal{op}}
\newcommand{\psiop}{\psi^{\op}}
\newcommand{\gammaop}{\gamma^{\op}}
\newcommand{\psiopforest}{F_{\psiop}}
\newcommand{\gammaopforest}{F_{\gammaop}}
\newcommand{\K}{K}

\newcommand{\FindStat}[1]{%
  \ifx&#1&%
  \url{www.findstat.org}
  \else \url{www.findstat.org/#1}
  \fi}%

\begin{document}

\title[Gorenstein dimensions of stable Auslander algebras]{Gorenstein dimensions for stable Auslander algebras of Nakayama algebras}
\date{}

\author[Giatagantzidis]{Odysseas Giatagantzidis}
\address[Giatagantzidis]{???}
\email{odysgiat@math.auth.gr}

\author[Hartung]{Carolin Hartung}
\address[Hartung]{Mathematical Institute of the University of Bonn, Germany}
\email{hartung@math.uni-bonn.de}

\author[Jain]{Nupur Jain}
\address[Jain]{Fakult\"at f\"ur Mathematik, Ruhr-Universit\"at Bochum, Germany}
\email{nupur.jain@rub.de}

\author[Marczinzik]{Ren\'e Marczinzik}
\address[Marczinzik]{Mathematical Institute of the University of Bonn, Germany}
\email{marczire@math.uni-bonn.de}

\author[Rubey]{Martin Rubey}
\address[Rubey]{Fakult\"at f\"ur Mathematik und Geoinformation, TU Wien, Austria}
\email{martin.rubey@tuwien.ac.at}

\author[Schmitz]{Sven Schmitz}
\address[Schmitz]{Mathematical Institute of the University of Bonn, Germany}
\email{sven.u.schmitz@uni-bonn.de}


\begin{abstract}
Let $A$ be a Nakayama algebra and $\sAus(A)$ its stable Auslander algebra. We introduce a new homological invariant, called the finitistic Gorenstein dimension, which coincides with the classical Gorenstein dimension whenever the latter is finite, and compute this invariant for $\sAus(A)$ when $A$ has infinite global dimension. Combined with an accompanying criterion for the finiteness of the Gorenstein dimension of $\sAus(A)$, our results provide a complete description of this dimension and verify the Gorenstein Symmetry Conjecture for this class of algebras. Furthermore, we show that the finitistic Gorenstein dimension of $ \sAus(A) $ is also left-right symmetric, is explicitly related to the finitistic dimension of $ A $, and can be attained among the injective dimensions of the projective $ \sAus(A) $-modules corresponding to simple $A$-modules. Our proofs rely on computing injective dimensions in terms of distances in the coresolution forest of $A$, a combinatorial structure closely related to Ringel's coresolution quiver.
\end{abstract}
\maketitle

\setcounter{tocdepth}{1}
\tableofcontents

\section{Introduction}
In their fundamental articles \cite{AR1,AR2,AR3,AR4,AR5,R1} on stable equivalences of dualizing $ R $-varieties, Auslander and Reiten established the basics on stable equivalences for module categories of Artin algebras, among other categories.
An important part of their theory concerns the homological dimensions of stable Auslander algebras in \cite{AR4} and the study of the stable module categories of Nakayama algebras in \cite{R1}; see also \cite{R2,R3} for further articles on stable module categories of Nakayama algebras.
Despite subsequent research on stable categories, it seems that no real progress has been made on the computation of the homological dimensions of stable Auslander algebras.

This paper builds on two foundational results from \cite{AR}.
In what follows, let $ A $ be a finite-dimensional algebra of finite representation type and let $\sAus(A):= \underline{\End}_A(X)$ be its \emph{stable Auslander algebra}, where $X$ is an additive generator of $\mod A$.
First, all indecomposable projective $ \sAus(A)$-modules were described along with injective coresolutions (see \Cref{theoreminjcoresstableaus}) that, crucially, are minimal when $ A $ is a Nakayama algebra. Furthermore, the homological complexity of $ \sAus(A) $ was bounded by that of $ A $ via the inequality $ \gldim \sAus(A) \leq 3 \gldim A - 1 $.

Our first main result associates the \emph{Gorenstein dimension} of $ \sAus(A) $, or rather a new finer homological invariant called the \Dfn{finitistic Gorenstein dimension}, with the finitistic dimension of $ A $, when $ A $ is a Nakayama algebra with infinite global dimension.
This new concept arises naturally in our context since it equals the maximum of the injective dimensions of projective $ \sAus(A) $-modules corresponding to simple $A$-modules, see \Cref{mainresult2}.
In general, we define the finitistic Gorenstein dimension of a ring $ R $ by
        \[
            \fingordim R = \sup \bigl\{ \idim P \mid P \in \proj R  \text{ with} \, \idim P < \infty \bigr\} ,
        \]
which should be interpreted as infinity when there is no projective $ R $-module with finite injective dimension. Recall that the Gorenstein dimension of $ R $, denoted by $ \gordim R $, is given by the above supremum waiving the requirement that $ \idim P $ be finite.

Specifically, our first main result is:

\begin{restatable}[{\Cref{cor:finGordim&findim,prop:findimA}}]{letterthm}{mainresultfour}\label{mainresult1}
Let $A$ be a Nakayama algebra with $ \gldim A = \infty $. Then:
        \[
            \fingordim \sAus(A)  =  \begin{cases}
                3 \findim A , & \text{if $ \, \findim A $ is even} \\
                3 \findim A + 1 , & \text{otherwise} ,
                    \end{cases}
        \]
        where the parity of $ \findim A $ is determined via a combinatorial criterion which also determines its exact value.
\end{restatable}

\Cref{mainresult1} yields a complete description of $ \gordim \sAus (A) $ when combined with \Cref{theo:three:components:left:vertex}, which provides a combinatorial characterization of the finiteness of $ \gordim \sAus(A) $ in this setting.

When $ A $ is additionally Iwanaga--Gorenstein, we obtain an even clearer picture:

\begin{lettercor}[\Cref{prop:iwanaga-gorenstein}]\label{lettercor}
    Let $ A $ be a Nakayama algebra with $ \gldim A = \infty $. If $A$ is Iwanaga--Gorenstein, then so is $\sAus(A)$, and
        \[
            \gordim \sAus(A) = 3\gordim A.
       \]
\end{lettercor}

This is a consequence of \Cref{mainresult1} combined with \cite[Remark 5.1.4]{Se}, but we also provide an alternative direct proof recovering the latter. Recall that a ring $ R $ is \emph{Iwanaga--Gorenstein} if both $ \gordim R $ and $ \gordim R^\op $ are finite, in which case they are equal \cite{Za}.

Another aspect of the finitistic Gorenstein dimension of $ \sAus(A) $ is that it is always attained on a certain subset of projective $ \sAus(A) $-modules:

\begin{restatable}[{\Cref{prop:infinite:height,theorem:simples,thm:idim-simple-finite}}]{letterprop}{mainresulttwo}    \label{mainresult2}
For every Nakayama algebra $A$,
    \[
        \fingordim \sAus(A) = \max \bigl\{ \idim \pstab{S} \mid S \in \mod A \textnormal{ is simple non-projective} \bigr\}   ,
    \]
where $ \pstab{M} := \underline{\Hom}_A(X,M) $ is the indecomposable projective $ \sAus(A) $-module corresponding to any indecomposable non-projective $ A $-module $ M $.
\end{restatable}

\Cref{mainresult2} does not extend to arbitrary algebras of finite representation type, as illustrated by \Cref{not.at.simples}.

The finitistic Gorenstein dimension of a ring $ R $ is inspired by the \emph{finitistic dimension} denoted by $ \findim R $, i.e., the supremum of the \emph{finite} projective dimensions of finitely generated $ R $-modules.
However, there is a fundamental difference at least when $ R $ is a semiperfect ring (e.g., a finite-dimensional algebra): $ \fingordim R = \infty $ means that all finitely generated projective $ R $-modules have infinite injective dimension, while $ \findim R = \infty $ means that, for every positive integer $ d $, there is a finitely generated $ R $-module with finite projective dimension exceeding $ d $. The former assertion holds since all projective $ R $-modules can be realized as direct sums of copies of the finitely many indecomposable projective $ R $-modules.

The finitistic dimension was first systematically studied by Hyman Bass in \cite{Ba}, where he also posed the question whether it is always finite. The question has a negative answer in general, as Nagata (see, for example, \cite{K}) and Kirkman--Kuzmanovich \cite{KK} showed with concrete counterexamples. Recent unverified computations announced by OpenAI claim to provide the first counterexamples in the class of finite-dimensional algebras in characteristic two and zero \cite{OAI1,OAI2}. These examples challenge the long-standing Finitistic Dimension Conjecture ($\mathsf{FDC}$), which asserts that the finitistic dimension is finite for every Artin algebra.

Another fundamental homological conjecture in the study of Artin algebras is the Gorenstein Symmetry Conjecture ($\mathsf{GSC}$) stating that $ \gordim A = \gordim { A^\op }$ for all Artin algebras $A$; see \cite[page 150]{AR6} and \cite[Conjecture 13]{ARS}.
The $\mathsf{GSC}$ is a consequence of the $\mathsf{FDC}$, and is still open since it is not affected by the above computational claims of OpenAI.

The second main result of the paper is the confirmation of the $\mathsf{GSC}$ for stable Auslander algebras of Nakayama algebras:

\begin{restatable}[\Cref{proposition:Gorenstein:symmetry}]{letterthm}{mainresultthree}
\label{mainresult3}
For every Nakayama algebra $A$,
\[
\gordim \sAus(A) = \gordim\sAus(A)^{\op}.
\]
\end{restatable}

Strangely enough, left-right symmetry also holds for the finitistic Gorenstein dimension of $ \sAus(A) $ (see \Cref{fingordim.symmetry}).
Observe that the two symmetry notions are not related a priori:
the algebra $ A $ of \Cref{fingordim.nonsymmetry}
has $ \gordim A = \infty = \gordim { A^\op } $ but $ \fingordim A \neq \fingordim { A^\op } $ (both finite).

Our methods entail the study of the \Dfn{coresolution forest} of a Nakayama algebra $ A $, which is an infinite periodic analog of the coresolution quiver of Gustafson, Madsen
and Ringel \cite{G,Ma,Ri1,Ri2}.
A fundamental property of this forest is that it can be embedded into the plane in an essentially unique way (see \Cref{thm:naturally}), allowing for the partition of its vertices into three disjoint subsets: vertebrae, left and right vertices. A further notion of importance is the maximal distance of a vertex from a vertebra, which we call the \Dfn{height of $ A $}, a quantity that goes back at least to \cite{G} and is closely related to the finitistic dimension of $ A $ (see \Cref{prop:findimA}).

We conclude this introduction by outlining the contents of the paper.
In \Cref{sec:preliminaries}, we review some preliminaries on homological dimensions, Nakayama algebras and general stable Auslander algebras.

In \Cref{sec:functional-digraph}, we introduce the \emph{coresolution forest} of a Nakayama algebra $ A $, a combinatorial object for computing injective dimensions of $ A $-modules, whose fundamental properties we settle in \Cref{fundamental.properties}. In \Cref{subsec:three-distance}, we establish the \emph{three-distance formula} which computes the injective dimensions of projective $ \sAus(A) $-modules by counting distances in the coresolution forest of $ A $. As a first application, we show that the injective dimension of a projective $ \sAus(A) $-module is either even or infinite. In \Cref{compute.fin.dim}, we compute the finitistic dimension of $ A $ in terms of its \emph{height}, an invariant of its coresolution forest.

In \Cref{sec:main-results-bounds}, we prove \Cref{mainresult1} and \Cref{lettercor} via a concise direct proof. Furthermore, we show that all indecomposable projective $ \sAus(A) $-modules corresponding to simple $ A $-modules have finite injective dimension (\Cref{thm:idim-simple-finite}), which yields half of \Cref{mainresult2} when combined with \Cref{mainresult1}, namely for $ \gldim A = \infty $. \Cref{sec:main-results-simples} is devoted to the proof of the second half of \Cref{mainresult2} for $ \gldim A < \infty $, which requires a more involved combinatorial treatment. In \Cref{sec:main-results-symmetry}, we establish a combinatorial criterion for the finiteness of $ \gordim \sAus(A) $, and then use it to prove \Cref{mainresult3}. We also prove the analogous left-right symmetry for the finitistic Gorenstein dimension of $ \sAus(A) $ (see \Cref{fingordim.symmetry}).

Finally, in \Cref{sec:further-questions}, we propose directions for future research based on our findings and provide examples illustrating the boundaries of potential generalizations.

\section{Preliminaries}\label{sec:preliminaries}
In this section, we collect useful facts that we use freely in the rest of the paper.

Let $A$ be a finite-dimensional algebra over a field $ \K $.
We denote by $ \mod A$ and $\proj A$ the category of finitely generated right $ A $-modules and finitely generated right projective $ A $-modules, respectively, and $D=\Hom_\K(-,\K)$ denotes the standard duality between finitely generated left and right $ A $-modules.
Unless stated otherwise, modules are finitely generated right modules.

\subsection{Homological dimensions}

Let $M$ be an $A$-module with minimal injective coresolution
$$0 \rightarrow M \rightarrow I^0 \rightarrow I^1 \rightarrow \cdots $$
and minimal projective resolution
$$\cdots \rightarrow P_1 \rightarrow P_0 \rightarrow M \rightarrow 0.$$
The \Dfn{injective dimension} $\idim M$ of $M$ is the smallest integer $n$, if it exists, such that $ I^{n+1} $ is zero.
The \Dfn{dominant dimension} $\domdim M$ of $M$ is the smallest integer $n$, if it exists, such that the module $I^n$ is not projective.
Both dimensions are infinite if no such $n$ exists.
In particular, the dominant dimension of a projective-injective module is infinite since zero modules are also projective-injective.
The \Dfn{projective dimension} $\pdim M$ and the \Dfn{codominant dimension} $\codomdim M$ of $M$ are defined dually.

The \Dfn{global dimension} of a finite-dimensional algebra $A$ is
\begin{align*}
    \gldim A &= \sup \bigl\{ \pdim M \mid M \in \mod A \bigr\}  \\
            &= \sup \bigl\{ \idim M \mid M \in \mod A \bigr\}   ,
\end{align*}
where the second equality is well-known, while its \Dfn{finitistic dimension} is
\[
  \findim A =\sup\bigl\{\pdim M \mid M\in\mod A \text{ with} \, \pdim M<\infty\bigr\}.
\]

Finally, the \Dfn{Gorenstein dimension} of $A$ is
\[\gordim A
= \idim A_A , \]
also called the self-injective dimension, for example, in \cite{Ba}.
It is immediate that $\gordim A \leq \gldim A$ for any finite-dimensional algebra $A$, and equality holds when $ \gldim A $ is finite; see, for example, \cite[Lemma~VI.5.5]{ARS}.

In this paper, we introduce a new homological dimension:

\begin{definition}
    The \Dfn{finitistic Gorenstein dimension} of a finite-dimensional algebra $ A $ is
        \[
            \fingordim A = \sup \bigl\{ \idim P \mid P \in \proj A  \text{ with} \, \idim P < \infty \bigr\}
        \] 
    which is interpreted as infinity when all projective $ A $-modules have infinite injective dimension.
\end{definition}

\begin{remark} \label{fingordim&findim}
    If $ \fingordim A $ is finite, then
        \[
            \fingordim A \leq \findim A^\op .
        \]
    Indeed, for any projective $ A $-module $ P $ with $ \idim P < \infty $, we have
        \[
            \idim_A P = \pdim_{A^\op} D(P) \leq \findim A^\op   .
        \]
    Moreover, equality holds when $ A $ is Iwanaga--Gorenstein as both dimensions are equal to the Gorenstein dimension, while \Cref{fGd.less.than.findim} illustrates that the inequality can also be proper.
\end{remark}

Note that the finitistic Gorenstein dimension is well-defined for any ring.

\subsection{Nakayama algebras}\label{subsec:preliminaries-nakayama}
This subsection collects useful well-known facts about Nakayama algebras.
For proofs and more details, we refer the reader to the textbooks \cite{AF,ARS,ASS,SkoYam}.

A \Dfn{Nakayama algebra} is a finite-dimensional $\K$-algebra $A$ such that every indecomposable $A$-module is \Dfn{uniserial}, meaning it has a unique composition series.

In this paper, all Nakayama algebras $A$ are connected and isomorphic to a bound quiver algebra.
Therefore, we assume that $A = \K Q/I$, with $Q$ being the linear
\begin{figure}[H]
\begin{center}
\begin{tikzpicture}[scale=1, >=stealth]
    \def\gap{1.2}

    \node (0) at (0,0) {$\mathbf{0}$};
    \node (1) at (\gap,0) {$\mathbf{1}$};
    \node (2) at (2*\gap,0) {$\mathbf{2}$};
    \node (n2) at (3.5*\gap,0) {$\mathbf{n-1}$};
    \node (n1) at (5.2*\gap,0) {$\mathbf{n}$};

    \draw[->, shorten <=1pt, shorten >=1pt] (0) -- (1);
    \draw[->, shorten <=1pt, shorten >=1pt] (1) -- (2);
    \draw[->, shorten <=1pt, shorten >=1pt] (n2) -- (n1);

    \draw[->, dotted, shorten <=1.5pt, shorten >=1.5pt] (2) -- (n2);

\end{tikzpicture}
\end{center}
\end{figure}
\noindent{}or cyclic
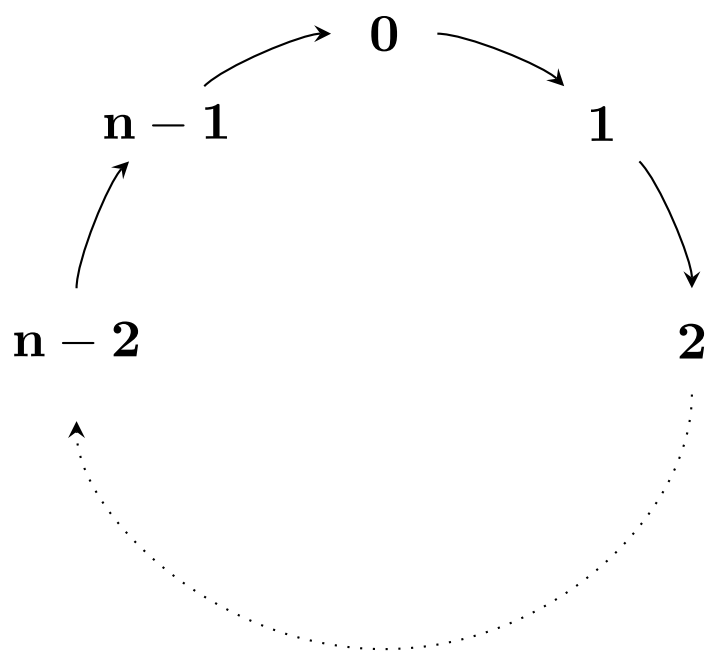
\begin{figure}[H]
\begin{center}
\begin{tikzpicture}[scale=0.55, >=stealth]
	\def\radius{3.7}

	\foreach \ang\lab in {
		90/0,
		45/1,
		0/2,
		315/{},
		270/{},
		225/{},
		180/{n-2},
		135/{n-1}
	}{
		\node (\lab) at (\ang:\radius) {$\mathbf{\lab}$};
	}
	\draw[->, shorten <=10pt, shorten >=10pt] (90:\radius) arc (90:45:\radius);    
	\draw[->, shorten <=10pt, shorten >=10pt] (45:\radius) arc (45:0:\radius);     
	\draw[->, dotted, shorten <=10pt, shorten >=15pt] (0:\radius) arc (360:180:\radius); 
    \draw[->, shorten <=10pt, shorten >=10pt] (180:\radius) arc (180:135:\radius);
    \draw[->, shorten <=10pt, shorten >=10pt] (135:\radius) arc (135:90:\radius);
\end{tikzpicture}
\end{center}
\end{figure}
\noindent{}quiver on $ n $ vertices $ 0 , 1 , \ldots , n-1 $. We call $A$ a \Dfn{linear Nakayama algebra} or a \Dfn{cyclic Nakayama algebra}, respectively. Observe that the arrows are of the form $ i \to i+1 $, except for $ n-1 \to 0 $ in the cyclic case. We also assume that $ A $ is not semisimple, i.e., we exclude the linear Nakayama algebras with a single vertex, since everything is well-understood for $ A \simeq \K $.

For each vertex $i$, let $ S_i $ be the simple $ A $-module supported at $ i $. Similarly, let $P_i = P(S_i) = e_i A$ and $I_i = I(S_i) = D(A e_i)$ denote the corresponding indecomposable projective and injective $A$-modules, respectively. A way to implement the simples is via the quotients $ S_i \cong e_i A / e_i J $, where $ J $ is the Jacobson radical of $ A $.

The module $P_i$ (respectively, $I_i$)  has a basis given by paths in $Q$ starting (respectively, ending) at $i$ that are not in the ideal $I$. Therefore, the vector space dimension of $P_i$ (respectively, $I_i$) is the length of the longest such path plus one, which we denote by $c_i$ (respectively, $d_i$). Consequently, $ c_i $ and $ d_i $ are the Loewy lengths of $ P_i $ and $ D(I_i) \cong A e_i $, respectively.
Collecting these numbers, we get the \Dfn{Kupisch} and \Dfn{coKupisch series} of $A$ denoted by $ c=\kupisch$ and $d=\cokupisch$, respectively. The entries of the two series satisfy $ c_{i} \leq c_{i + 1} + 1 $ and $d_{i} + 1 \geq d_{i + 1}$ for all $ i $, where $ c_i \geq 2 $ (respectively, $ d_i \geq 2 $) for all $ i $ except for $ i=n-1 $ (respectively, $ i= 0 $) in the linear case; then $ c_{n-1} = 1 $ and $ d_0 = 1 $.

Every sequence of positive integers satisfying one of these two sets of inequalities determines a unique connected Nakayama algebra; see \cite[Theorem 32.9]{AF}.
Furthermore, it determines a unique sequence of the other kind according to the next lemma (and its dual):
\begin{lemma}[\mbox{\cite[Theorem 2.2]{F}}]\label{lem:kupisch-cokupisch}
    The Kupisch series and the coKupisch series of $A$ are related via
    \[
    c_i = \max \bigl\{ k \in \Z \mid k\leq d_{i+k-1} \bigr\}.
    \]
\end{lemma}
\begin{proof}
    Observe that $c_i-1$ is the maximal length of a nonzero path in $Q$ starting at vertex $ i $. Therefore, for a positive integer $ k $, $ k \leq c_i $ if and only if the path of length $ k-1 $ starting at vertex $ i $ is nonzero. Dually, this is equivalent to $ k \leq d_{ i+k-1 } $.
\end{proof}

\begin{example}\label{example:linearnakayama}
    Let $Q$ be the quiver
    \begin{figure}[H]
    \begin{center}
        \begin{tikzpicture}[>=stealth]
            \def\gap{1.2}
            \node (0) at (0*\gap,0) {$\mathbf{0}$};
            \node (1) at (1*\gap,0) {$\mathbf{1}$};
            \node (2) at (2*\gap,0) {$\mathbf{2}$};
            \node (3) at (3*\gap,0) {$\mathbf{3}$};
            \node (4) at (4*\gap,0) {$\mathbf{4}$};
            \node (5) at (5*\gap,0) {$\mathbf{5}$};
            \node (6) at (6*\gap,0) {$\mathbf{6}$};
            \draw[->, shorten <=1pt, shorten >=1pt] (0) -- node[above] {$a$} (1);
            \draw[->, shorten <=1pt, shorten >=1pt] (1) -- node[above] {$b$} (2);
            \draw[->, shorten <=1pt, shorten >=1pt] (2) -- node[above] {$c$} (3);
            \draw[->, shorten <=1pt, shorten >=1pt] (3) -- node[above] {$d$} (4);
            \draw[->, shorten <=1pt, shorten >=1pt] (4) -- node[above] {$e$} (5);
            \draw[->, shorten <=1pt, shorten >=1pt] (5) -- node[above] {$f$} (6);
        \end{tikzpicture}
    \end{center}
    \end{figure}
    \noindent{}bound by the admissible ideal $I = \langle abcd, cdef\rangle$. Let $A = \K Q/I$ be the respective linear Nakayama algebra.
    Counting the length of the longest nonvanishing path in $ A $ starting at every vertex $i$, we get the Kupisch series $c = [4,5,4,4,3,2,1]$. Similarly, one computes the coKupisch series $d = [1,2,3,4,4,5,4]$ either by counting the length of the longest nonvanishing path in $A$ ending at every $i$, or by the dual of \Cref{lem:kupisch-cokupisch}.
    \end{example}

Every indecomposable $A$-module $M$ for a Nakayama algebra $A$ is isomorphic to $ \indec $ for unique integers $ i $ and $ k $ such that $0 \leq i \leq n - 1$ and $1 \leq k \leq c_i$. Let $\mik := \indec$ denote this module, i.e., the unique indecomposable $A$-module with top isomorphic to $S_i$ and Loewy length $k$.

We close this subsection with a description of the projective cover $P(M)$ and the injective envelope $I(M)$ of an indecomposable $ A $-module $ M $, as well as of its syzygy $\Omega(M)$ and cosyzygy $\Sigma (M)$.

\begin{proposition}\label{prop:injenvcosyz}
Let $A$ be a Nakayama algebra and  $M = \mik$ an indecomposable $A$-module. Then:
\begin{enumerate}
    \item $P(M) = \mik[i][c_i] $ and $\Omega(M) = \mik[i + k][c_i - k]$.  In particular, $M(i,k)$ is projective if and only if $k=c_i$.
    \item $ I(M) = \mik[i + k - d_{i + k - 1}][d_{i + k - 1}] $ and $\cosyz M = \mik[i + k - d_{i + k - 1}][d_{i + k - 1} - k]$. In particular, $M(i,k)$ is injective if and only if $k=d_{i + k - 1}$.
\end{enumerate}
\end{proposition}
The proof is an elementary translation of the homological terms into combinatorial terms. We refer, for example, to \cite{MRS} for details.

\subsection{Homological algebra of stable Auslander algebras}
In this subsection we collect some results from the survey article \cite{AR}. We will restrict to representation-finite algebras for simplicity. For a more general functorial approach, we refer to \cite{AR,AR1} and, for a textbook introduction, to \cite[Section~5.11]{Zi}.

Let $A$ be a finite-dimensional algebra.
Given two $A$-modules $M$ and $N$, we denote by $\underline{\Hom}_A(M,N)$ the $\K$-vector space of $A$-linear maps from $M$ to $N$, modulo the subspace consisting of maps that factor through a projective module. If $ A $ is representation-finite, then the stable Auslander algebra of $A$ is defined as $\sAus(A):=\underline{\End}_A(X)$, where $X$ is the direct sum of all indecomposable $A$-modules.

\begin{theorem} \label{theoreminjcoresstableaus}
Let $ A $ be a finite-dimensional algebra of finite representation type and $\sAus(A) = \underline{\End}_A(X)$ its stable Auslander algebra. Let $\pstab M := \underline{\Hom}_A(X,M)$ denote the indecomposable projective $\sAus(A)$-module induced by the indecomposable non-projective $A$-module $M$. Then:
\begin{enumerate}
    \item \label{theoreminjcoresstableaus1}Every indecomposable projective $\sAus(A)$-module is isomorphic to $\pstab M$ for a unique $M$.
    \item \label{theoreminjcoresstableaus2}$\pstab M$ is injective if and only if $M$ has codominant dimension at least one.
    \item \label{theoreminjcoresstableaus3}There is an injective coresolution of $\pstab M$ of the form:
    \begin{align*}
        0 \rightarrow \underline{\Hom}_A(X,M) \rightarrow \Ext_A^1(X,\Omega(M)) \rightarrow \Ext_A^1(X,P(M)) \rightarrow \Ext_A^1(X,M)\\
        \rightarrow \Ext_A^2(X,\Omega(M)) \rightarrow \Ext_A^2(X,P(M)) \rightarrow \Ext_A^2(X,M) \rightarrow \cdots.
    \end{align*}

\end{enumerate}
\end{theorem}
\begin{proof}
For (1), see \cite[Proposition I.3.3]{AR}.
For (2), see \cite[Proposition II.1.2]{AR} and for (3) see the proof of \cite[Proposition I.3.2]{AR}.
\end{proof}

\begin{remark} \label{injresremark}
    The following dual version of \Cref{theoreminjcoresstableaus} also holds:
    \begin{enumerate}
        \item\label{rem:injective.sAus} Every indecomposable injective $\sAus(A)$-module is isomorphic to $\Ext_A^1(X,M)$ for a unique indecomposable non-injective $A$-module $M$.
        \item $\Ext_A^1(X,M)$ is projective if and only if $M$ has dominant dimension at least one.
        \item There is a projective resolution of $ \Ext_A^1(X,M)$ of the form:
        \begin{multline*}
            \cdots \to \underline{\Hom}_A(X,\Omega^{1}(M))
            \to \underline{\Hom}_A(X,\Omega^{1}(I(M)))  \\
            \to \underline{\Hom}_A(X,\Omega \Sigma(M))
            \to \underline{\Hom}_A(X,M) \\
            \to \underline{\Hom}_A(X,I(M))
            \to \underline{\Hom}_A(X,\Sigma^1(M))
            \to \Ext_A^1(X,M)
            \to 0.
        \end{multline*}
    \end{enumerate}
\end{remark}

We conclude this section with two useful facts.
First, the stable Auslander algebra of a disconnected algebra is isomorphic to the direct product of the stable Auslander algebras of its blocks, so restricting to connected Nakayama algebras incurs no loss of generality.

The second fact we record as a proposition as it will be used in \Cref{sec:main-results-symmetry}:

\begin{proposition}\label{Aus.op}
    For every finite-dimensional algebra $ A $ of finite representation type,
        \[
            \sAus(A)^{\op} \cong \sAus(A^{\op}).
        \]
\end{proposition}

\begin{proof}
    Direct consequence of the transpose duality $ \underline{\mod} A \to \underline{ \mod } A^\op $ as the transpose of $ X $ is an additive generator of $ \underline{\mod} A^\op $; see \cite[Proposition IV.1.6.(b)]{ARS}.
\end{proof}

\section{Injective dimensions via coresolution forests}\label{sec:coresolution-forests}

In this section, we investigate the combinatorial properties of the coresolution forest of a Nakayama algebra and use it to compute homological dimensions over the algebra and its
stable Auslander algebra. This forest is a variation of the
coresolution quiver introduced by Madsen and Ringel~\cite{Ma,Ri1,Ri2}.

\subsection{Coresolution forests of \texorpdfstring{$n$}{n}-coKupisch series}\label{sec:functional-digraph}

We begin with an extension of Kupisch and coKupisch series to periodic doubly infinite sequences:

\begin{definition}
  A sequence $(c_i)_{i\in\Z}$ of positive integers is an
  \Dfn{$n$-Kupisch series} if $c_i\leq c_{i+1} + 1$ and
  $c_i = c_{i+n}$.  A sequence $(d_i)_{i\in\Z}$ of positive integers
  is an \Dfn{$n$-coKupisch series} if $d_{i+1}\leq d_i + 1$ and
  $d_i = d_{i+n}$. The \Dfn{associated $ n $-Kupisch series} of an $ n $-coKupisch series is the one determined uniquely by \Cref{lem:kupisch-cokupisch}.
\end{definition}

Every $ n $-coKupisch series comes together with a function $ \psi \colon \Z \to \Z $ and two associated functional digraphs as follows:

\begin{definition}\label{dfn:coresolution-forest}
    The \Dfn{coresolution function} of an $ n $-coKupisch series $(d_i)_{i\in\Z}$ is the function $ \psi \colon \Z \to \Z $ with $ \psi(i) = i - d_i $ for all $ i \in \Z $. Its \Dfn{coresolution forest} $ \psiforest $ is the digraph with vertex set $\Z$ and edges
  $i \rightarrow \psi(i)$ for all $i \in \Z$.  Finally, its \Dfn{coresolution quiver}, denoted by $ \psiquiver $, is the functional digraph of $ \psi $ viewed as a function $ \Z_n \to \Z_n $.
\end{definition}

      Note that $\psiforest$ is indeed an (oriented) forest since every vertex has a unique parent and $\psi$ is strictly decreasing; see \Cref{prop:cokupisch}.

Observe that our definition of $ \psi $ does not agree with the established definition for Nakayama algebras in the literature, but rather coincides with the dual function $ \gamma $. This is because we deal with right modules by default, contrary to the prevalent convention in papers such as \cite{Ma,Ri2,S,S2}. In particular, we are justified in calling $ \psi $ the \emph{coresolution} function and $ Q_\psi $ the \emph{coresolution} quiver of $(d_i)_{i\in\Z}$, since $ \psi $ plays the same role for right modules as $ \gamma $ does for left modules.

We define now this dual function $ \gamma $ and the related digraphs as we will need them later. A finite version of the resolution forest for linear Nakayama algebras was studied in \cite{KMMRS}.

\begin{definition}\label{dfn:resolution-forest}
    The \Dfn{resolution function} of an $ n $-coKupisch series $(d_i)_{i\in\Z}$ is the function $ \gamma \colon \Z \to \Z $ defined by $ \gamma(i) = i + c_i $, where $(c_i)_{i\in\Z}$ is the $ n $-Kupisch series associated to $(d_i)_{i\in\Z}$. Furthermore, the \Dfn{resolution forest} $\gammaforest $ and \Dfn{resolution quiver} $ Q_\gamma $ of $(d_i)_{i\in\Z}$ are the functional digraphs of $ \gamma $ and of its induced function $ \Z_n \to \Z_n $, respectively.
\end{definition}

Next, we provide a first indication of the close connection between an $ n $-coKupisch series and its coresolution function, whose proof is elementary and thus omitted:

\begin{lemma}\label{prop:cokupisch}
  A sequence $(d_i)_{i\in\Z}$ is an $n$-coKupisch series if and only
  if its coresolution function $\psi:\Z\to\Z$ defined by $\psi(i) := i-d_i$ is:
  \begin{enumerate}
  \item\label{prop:order-preserving} order preserving, that is, $i\leq j\Rightarrow\psi(i)\leq\psi(j)$;
  \item strictly decreasing, that is, $\psi(i)<i$; and
  \item $i - \psi(i)$ is $n$-periodic, that is, $\psi(n+i)-\psi(i)=n$.
  \end{enumerate}
\end{lemma}

Throughout the paper, we use standard terminology of forests. Specifically, the \Dfn{parent} of a vertex $i \in \psiforest$ is $\psi(i)$, and $i$ is a \Dfn{child} of $\psi(i)$. A vertex without children is a \Dfn{leaf}. A vertex $u$ is an \Dfn{ancestor} of $v$ if there is a directed (possibly trivial) path from $v$ to $u$, i.e., if $\psi^\ell(v) = u$ for some integer $\ell \geq 0$.
If $u$ is an ancestor of $v$ and $u \neq v$, we call $u$ a \Dfn{proper ancestor} of $v$.
Finally, a \Dfn{component} of $\psiforest$ is a maximal connected subgraph, and every component is a (directed) \Dfn{tree}.

In what follows, we draw parents above their children in a coresolution forest, omitting arrowheads since edges are implicitly directed from bottom to top. Furthermore, we indicate negative integers with an overline, so that $-1$ is written as $\overline{1}$ and so on.

\begin{example}\label{ex:psiforest}
    In \Cref{fig:coresforest-1234454}, we show a portion of the coresolution forest of the linear Nakayama algebra discussed in \Cref{example:linearnakayama}. The colored vertices illustrate the computations carried out in \Cref{ex:3.dist.form}.
    \Cref{fig:coresforest-55667} shows a portion of the coresolution forest of the cyclic Nakayama algebra with coKupisch series $d = [5,5,6,6,7]$.
    \begin{figure}[t]
      \begin{subfigure}[t]{0.4\textwidth}
      \centering
  \centering
  {
  \newcommand{\nodea}{\node[draw=none, inner sep=2pt] (a) {}
  ;}
  \newcommand{\nodeb}{\node[draw=none, inner sep=2pt] (b) {}
  ;}
  \newcommand{\nodec}{\node[draw=none, inner sep=2pt] (c) {}
  ;}
  \newcommand{\noded}{\node[draw=none, inner sep=2pt] (d) {}
  ;}
  \newcommand{\nodee}{\node (e) {$\overline{10}$}
  ;}
  \newcommand{\nodef}{\node (f) {$\overline{9}$}
  ;}
  \newcommand{\nodeg}{\node (g) {$\overline{8}$}
  ;}
  \newcommand{\nodeh}{\node (h) {$\overline{7}$}
  ;}
  \newcommand{\nodei}{\node (i) {$\overline{6}$}
  ;}
  \newcommand{\nodej}{\node (j) {$\overline{5}$}
  ;}
  \newcommand{\nodek}{\node (k) {$\overline{4}$}
  ;}
  \newcommand{\nodel}{\node (l) {$\overline{3}$}
  ;}
  \newcommand{\nodem}{\node (m) {$\overline{2}$}
  ;}
  \newcommand{\noden}{\node (n) {$\overline{1}$}
  ;}
  \newcommand{\nodeo}{\node (o) {$0$}
  ;}
  \newcommand{\nodep}{\node (p) {$1$}
  ;}
  \newcommand{\nodeq}{\node (q) {$2$}
  ;}
  \newcommand{\noder}{\node (r) {$3$}
  ;}
  \newcommand{\nodes}{\node (s) {$4$}
  ;}
  \newcommand{\nodet}{\node (t) {$5$}
  ;}
  \newcommand{\nodeu}{\node (u) {$6$}
  ;}
  \newcommand{\nodev}{\node (v) {$7$}
  ;}
  \newcommand{\nodew}{\node (w) {$8$}
  ;}
  \newcommand{\nodex}{\node (x) {$9$}
  ;}
  \newcommand{\nodey}{\node (y) {$10$}
  ;}
  \newcommand{\nodez}{\node (z) {$11$}
  ;}
  \newcommand{\nodeaa}{\node (aa) {$12$}
  ;}
  \newcommand{\nodeab}{\node (ab) {$13$}
  ;}
  \newcommand{\nodeac}{\node[draw=none, inner sep=2pt] (ac) {}
  ;}
  \newcommand{\nodead}{\node[draw=none, inner sep=2pt] (ad) {}
  ;}
  \newcommand{\nodeae}{\node[draw=none, inner sep=2pt] (ae) {}
  ;}
  \newcommand{\nodeaf}{\node[draw=none, inner sep=2pt] (af) {}
  ;}
  \begin{tikzpicture}[auto, font=\footnotesize]
  \matrix[column sep=.3cm, row sep=.5cm,ampersand replacement=\&]{
    \nodea  \&  \nodeb  \&  \nodec  \&  \noded  \\
    \nodee  \&  \nodef  \&  \nodeg  \&           \\
    \nodeh  \&  \nodei  \&  \nodej  \&  \nodek  \\
    \nodel  \&  \nodem  \&  \noden  \&           \\
    \nodeo  \&  \nodep  \&  \nodeq  \&  \noder  \\
    \nodes  \&  \nodet  \&  \nodeu  \&           \\
    \nodev  \&  \nodew  \&  \nodex  \&  \nodey  \\
    \nodez  \&  \nodeaa  \&  \nodeab  \&           \\
    \nodeac  \&  \nodead  \&  \nodeae  \&  \nodeaf  \\
  };

  \path[ultra thick]   (g) edge (h)
  	  (g) edge (k)
  	  (h) edge (l)
  	  (h) edge (m)
  	  (j) edge (n)
  	  (n) edge (o)
  	  (n) edge (p)
  	  (n) edge (q)
  	  (n) edge (r)
  	  (o) edge (s)
  	  (o) edge (t)
  	  (u) edge (v)
  	  (u) edge (w)
  	  (u) edge (x)
  	  (u) edge (y)
  	  (x) edge (ab);
  \path[ultra thick, dotted]   (a) edge (e)
  	  (a) edge (f)
  	  (c) edge (g)
  	  (ab) edge (ac)
  	  (ab) edge (ad)
  	  (ab) edge (ae)
  	  (ab) edge (af);
  \path[ultra thick, draw=green]  (v) edge (z)  (v) edge (aa);
  \path[ultra thick, draw=red]    (g) edge (i)  (g) edge (j);
  \path[ultra thick, draw=blue]   (q) edge (u);
  \begin{pgfonlayer}{background}
  \foreach \v in {i,j,m}{\node[circle, fill=red, opacity=0.35, inner sep=4pt] at (\v) {};}
  \foreach \v in {q,r,u}{\node[circle, fill=blue, opacity=0.35, inner sep=4pt] at (\v) {};}
  \foreach \v in {z,aa,w}{\node[circle, fill=green, opacity=0.35, inner sep=4pt] at (\v) {};}
  \end{pgfonlayer}
  \end{tikzpicture}}
        \caption{$d = [1,2,3,4,4,5,4]$}
        \label{fig:coresforest-1234454}
      \end{subfigure}\hfill
      \begin{subfigure}[t]{0.4\textwidth}
      \centering
        \input{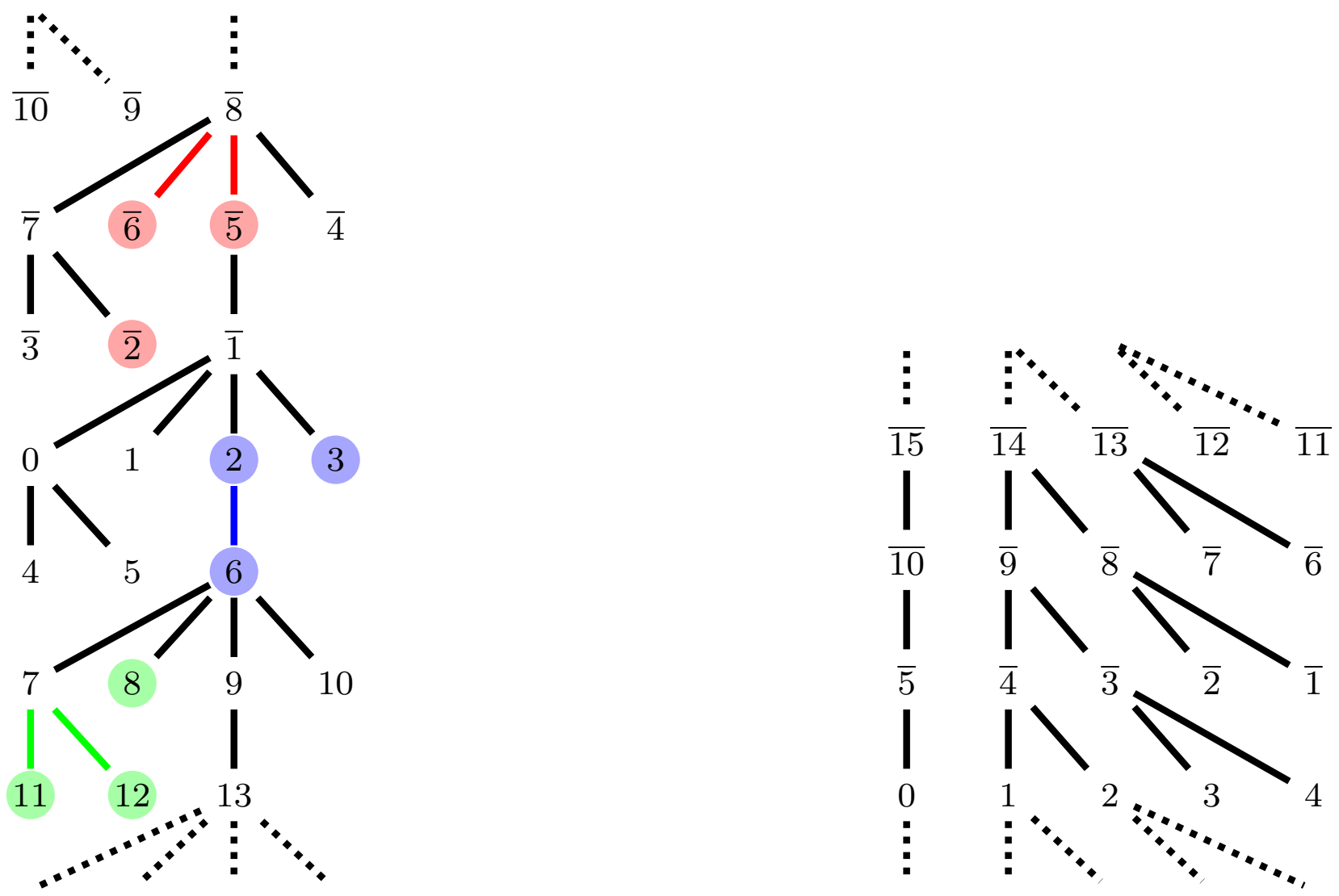}
        \caption{$d = [5, 5, 6, 6, 7]$}
        \label{fig:coresforest-55667}
      \end{subfigure}\hfill
      \caption{Portions of coresolution forests of the two Nakayama algebras with the displayed coKupisch series.}
      \label{fig:cores-forests-examples}
\end{figure}
\end{example}

  Let $\rho \colon \psiforest \to \psiquiver $ be the map sending each vertex of the coresolution forest to
  its residue class modulo $n$.
  It is straightforward to check that $\rho$ is a covering map since $d_i = i-\psi(i)$ is $n$-periodic.
  That is, $\rho$ is
  surjective and maps the parent and children of a vertex $i \in \psiforest$ bijectively to the successor and the predecessors of
  $\rho(i) \in \psiquiver$, respectively. More precisely, each component of the coresolution forest is a universal
  cover of the component in the coresolution quiver to which it maps.

\begin{proposition}\label{prop:covering}
The preimage of every oriented cycle of $\psiquiver$ under $\rho$ is a collection of doubly infinite oriented paths. Moreover, each component of $\psiforest$ contains a unique such path.
\end{proposition}

\begin{proof}
The result follows from the preceding discussion and the observation that every vertex on an oriented cycle of $\psiquiver$ has a unique predecessor that lies on that cycle.
\end{proof}

\begin{definition}
    The unique doubly infinite oriented path in every component of $ \psiforest $ is called its \Dfn{spine}. A vertex lying on the spine is a \Dfn{vertebra}.
\end{definition}

Following the terminology of \cite[Section 3]{Ma}, for a Nakayama algebra $ A $, the vertebrae of $ \psiforest $ are exactly the vertices corresponding to $ \psi $-regular simple $ A $-modules. It is also proven there that these are exactly the simple $ A $-modules with even projective dimension when $ \gldim A < \infty $, or else with infinite projective dimension.

  The number of components in $ \psiforest $ mapping onto a component $ C $ of $ \psiquiver $ is
  \[
    \frac{1}{n}\sum_{j=1}^\ell d_{i_j},
  \]
  where $i_1,\dots,i_\ell$ are the vertices of $\psiquiver$ in the (unique) cycle of $ C $.  This quantity is called the weight of the component $
  C $, and $ \ell $ is its length.  It turns out that neither quantity depends on $ C $, see \cite[Proposition 1.1]{S}.  Therefore we define the \Dfn{weight} and the \Dfn{length} of $(d_i)_{i\in\Z}$ (or of the Nakayama algebra $ A $ with this $ n $-coKupisch series) as the weight and the length of any component of $\psiquiver$, respectively. In particular, we have:

  \begin{lemma}\label{rem:weight}
      The number of components in $ \psiforest $ is equal to the number of components in $ \psiquiver $ multiplied by the weight of $(d_i)_{i\in\Z}$.
  \end{lemma}

\begin{example}\label{ex:fig.2}
  \Cref{fig:coresforest-8999} shows the coresolution quiver and
  coresolution forest of the Nakayama algebra $ A $ with coKupisch series
  $d = [8,9,9,9]$.  The quiver is connected, and the weight of its
  cycle is two. Therefore, the coresolution forest has two (isomorphic) components, both
  of which are universal covers of the quiver.
    \begin{figure}[t]
      \centering
        \begin{minipage}[c]{0.4\textwidth}
    \centering
    \begin{tikzpicture}[>=stealth, font=\footnotesize,
        node distance=1.2cm and 1.4cm]
      \node (q0) {$0$};
      \node (q1) [right of=q0] {$1$};
      \node (q2) [right of=q1] {$2$};
      \node (q3) [right of=q2] {$3$};
      \draw[->] (q1) -- (q0);
      \draw[->] (q2) to (q1);
      \draw[->] (q3) to (q2);
      \draw[->] ([yshift=5.2pt]q0) arc[start angle=25, end angle=335, radius=4.5mm];
    \end{tikzpicture}
  \end{minipage}%
  \hfill
  \begin{minipage}[c]{0.55\textwidth}
\centering
{
\newcommand{\nodea}{\node[draw=none, inner sep=2pt] (a) {}
;}
\newcommand{\nodeb}{\node[draw=none, inner sep=2pt] (b) {}
;}
\newcommand{\nodec}{\node[draw=none, inner sep=2pt] (c) {}
;}
\newcommand{\noded}{\node[draw=none, inner sep=2pt] (d) {}
;}
\newcommand{\nodee}{\node[draw=none, inner sep=2pt] (e) {}
;}
\newcommand{\nodef}{\node[draw=none, inner sep=2pt] (f) {}
;}
\newcommand{\nodeg}{\node[draw=none, inner sep=2pt] (g) {}
;}
\newcommand{\nodeh}{\node[draw=none, inner sep=2pt] (h) {}
;}
\newcommand{\nodei}{\node (i) {$\overline{16}$}
;}
\newcommand{\nodej}{\node (j) {$\overline{15}$}
;}
\newcommand{\nodek}{\node (k) {$\overline{14}$}
;}
\newcommand{\nodel}{\node (l) {$\overline{13}$}
;}
\newcommand{\nodem}{\node (m) {$\overline{12}$}
;}
\newcommand{\noden}{\node (n) {$\overline{11}$}
;}
\newcommand{\nodeo}{\node (o) {$\overline{10}$}
;}
\newcommand{\nodep}{\node (p) {$\overline{9}$}
;}
\newcommand{\nodeq}{\node (q) {$\overline{8}$}
;}
\newcommand{\noder}{\node (r) {$\overline{7}$}
;}
\newcommand{\nodes}{\node (s) {$\overline{6}$}
;}
\newcommand{\nodet}{\node (t) {$\overline{5}$}
;}
\newcommand{\nodeu}{\node (u) {$\overline{4}$}
;}
\newcommand{\nodev}{\node (v) {$\overline{3}$}
;}
\newcommand{\nodew}{\node (w) {$\overline{2}$}
;}
\newcommand{\nodex}{\node (x) {$\overline{1}$}
;}
\newcommand{\nodey}{\node (y) {$0$}
;}
\newcommand{\nodez}{\node (z) {$1$}
;}
\newcommand{\nodeaa}{\node (aa) {$2$}
;}
\newcommand{\nodeab}{\node (ab) {$3$}
;}
\newcommand{\nodeac}{\node (ac) {$4$}
;}
\newcommand{\nodead}{\node (ad) {$5$}
;}
\newcommand{\nodeae}{\node (ae) {$6$}
;}
\newcommand{\nodeaf}{\node (af) {$7$}
;}
\newcommand{\nodeag}{\node (ag) {$8$}
;}
\newcommand{\nodeah}{\node (ah) {$9$}
;}
\newcommand{\nodeai}{\node (ai) {$10$}
;}
\newcommand{\nodeaj}{\node (aj) {$11$}
;}
\newcommand{\nodeak}{\node (ak) {$12$}
;}
\newcommand{\nodeal}{\node (al) {$13$}
;}
\newcommand{\nodeam}{\node (am) {$14$}
;}
\newcommand{\nodean}{\node (an) {$15$}
;}
\newcommand{\nodeao}{\node (ao) {$16$}
;}
\newcommand{\nodeap}{\node (ap) {$17$}
;}
\newcommand{\nodeaq}{\node (aq) {$18$}
;}
\newcommand{\nodear}{\node (ar) {$19$}
;}
\newcommand{\nodeas}{\node (as) {$20$}
;}
\newcommand{\nodeat}{\node (at) {$21$}
;}
\newcommand{\nodeau}{\node (au) {$22$}
;}
\newcommand{\nodeav}{\node (av) {$23$}
;}
\newcommand{\nodeaw}{\node[draw=none, inner sep=2pt] (aw) {}
;}
\newcommand{\nodeax}{\node[draw=none, inner sep=2pt] (ax) {}
;}
\newcommand{\nodeay}{\node[draw=none, inner sep=2pt] (ay) {}
;}
\newcommand{\nodeaz}{\node[draw=none, inner sep=2pt] (az) {}
;}
\newcommand{\nodeba}{\node[draw=none, inner sep=2pt] (ba) {}
;}
\newcommand{\nodebb}{\node[draw=none, inner sep=2pt] (bb) {}
;}
\newcommand{\nodebc}{\node[draw=none, inner sep=2pt] (bc) {}
;}
\newcommand{\nodebd}{\node[draw=none, inner sep=2pt] (bd) {}
;}
\begin{tikzpicture}[auto, font=\footnotesize]
\matrix[column sep=.3cm, row sep=.5cm,ampersand replacement=\&]{
  \nodea  \&  \nodeb  \&  \nodec  \&  \noded  \&  \nodee  \&  \nodef  \&  \nodeg  \&  \nodeh  \\
  \nodei  \&  \nodej  \&  \nodek  \&  \nodel  \&  \nodem  \&  \noden  \&  \nodeo  \&  \nodep  \\
  \nodeq  \&  \noder  \&  \nodes  \&  \nodet  \&  \nodeu  \&  \nodev  \&  \nodew  \&  \nodex  \\
  \nodey  \&  \nodez  \&  \nodeaa  \&  \nodeab  \&  \nodeac  \&  \nodead  \&  \nodeae  \&  \nodeaf  \\
  \nodeag  \&  \nodeah  \&  \nodeai  \&  \nodeaj  \&  \nodeak  \&  \nodeal  \&  \nodeam  \&  \nodean  \\
  \nodeao  \&  \nodeap  \&  \nodeaq  \&  \nodear  \&  \nodeas  \&  \nodeat  \&  \nodeau  \&  \nodeav  \\
  \nodeaw  \&  \nodeax  \&  \nodeay  \&  \nodeaz  \&  \nodeba  \&  \nodebb  \&  \nodebc  \&  \nodebd  \\
};

\path[ultra thick]   (i) edge (q)
	  (i) edge (r)
	  (j) edge (s)
	  (k) edge (t)
	  (m) edge (u)
	  (m) edge (v)
	  (n) edge (w)
	  (o) edge (x)
	  (q) edge (y)
	  (q) edge (z)
	  (r) edge (aa)
	  (s) edge (ab)
	  (u) edge (ac)
	  (u) edge (ad)
	  (v) edge (ae)
	  (w) edge (af)
	  (y) edge (ag)
	  (y) edge (ah)
	  (z) edge (ai)
	  (aa) edge (aj)
	  (ac) edge (ak)
	  (ac) edge (al)
	  (ad) edge (am)
	  (ae) edge (an)
	  (ag) edge (ao)
	  (ag) edge (ap)
	  (ah) edge (aq)
	  (ai) edge (ar)
	  (ak) edge (as)
	  (ak) edge (at)
	  (al) edge (au)
	  (am) edge (av);
\path[ultra thick, dotted]   (a) edge (i)
	  (a) edge (j)
	  (b) edge (k)
	  (c) edge (l)
	  (e) edge (m)
	  (e) edge (n)
	  (f) edge (o)
	  (g) edge (p)
	  (ao) edge (aw)
	  (ao) edge (ax)
	  (ap) edge (ay)
	  (aq) edge (az)
	  (as) edge (ba)
	  (as) edge (bb)
	  (at) edge (bc)
	  (au) edge (bd);
\path[ultra thick, draw=red]    (i) edge (r)  (r) edge (aa)  (aa) edge (aj);
\path[ultra thick, draw=green]  (m) edge (v)  (v) edge (ae)  (ae) edge (an);
\path[ultra thick, draw=blue]   (i) edge (q)  (q) edge (z)  (z) edge (ai)  (ai) edge (ar);
\begin{pgfonlayer}{background}
\node[circle, fill=red, opacity=0.35, inner sep=4pt] at (aj) {};
\node[circle, fill=green, opacity=0.35, inner sep=4pt] at (an) {};
\node[circle, fill=blue, opacity=0.35, inner sep=4pt] at (ar) {};
\end{pgfonlayer}
\end{tikzpicture}}
  \end{minipage}%
      \caption{The coresolution quiver (left) and the coresolution forest (right) of the Nakayama algebra $ A $ with coKupisch series $ d = [8,9,9,9]$. The length of $ A $ is one, and its weight is $ 8 / 4=2 $.}
      \label{fig:coresforest-8999}
    \end{figure}
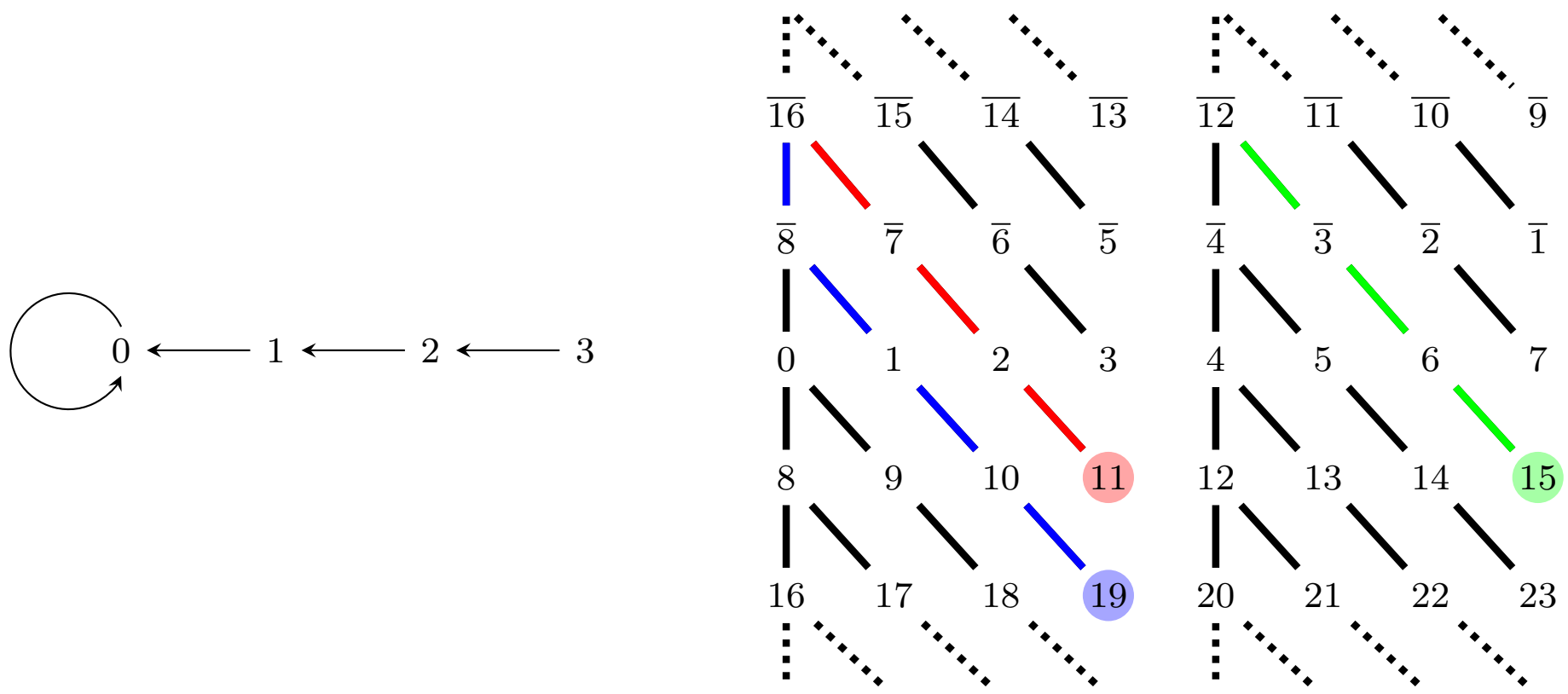
\end{example}

\begin{example}
  Similarly to \Cref{ex:fig.2}, \Cref{fig:coresforest-5445455654} shows the coresolution quiver $ \psiquiver $ and the coresolution forest $ \psiforest $ of the Nakayama algebra $ A $ with coKupisch series $d = [5,4,4,5,4,5,5, \allowbreak 6,5,4]$.
    \begin{figure}[t]
        \centering
          \centering
  \begin{minipage}[c]{0.4\textwidth}
    \centering
    \begin{tikzpicture}[>=stealth, font=\footnotesize,
        every node/.style={inner sep=1.5pt}]
      \node (n4) at (0,0)     {$4$};
      \node (n0) at (1.2,0)   {$0$};
      \node (n5) at (2.4,0)   {$5$};
      \node (n9) at (3.6,0)   {$9$};
      \node (n6) at (0,-1.4)   {$6$};
      \node (n1) at (1.2,-1.4) {$1$};
      \node (n7) at (2.4,-1.4) {$7$};
      \node (n2) at (0,-2.8)   {$2$};
      \node (n8) at (1.2,-2.8) {$8$};
      \node (n3) at (2.4,-2.8) {$3$};

      \draw[->] (n4) -- (n0);
      \draw[->] (n9) -- (n5);
      \draw[->] (n0) to[bend left=45] (n5);
      \draw[->] (n5) to[bend left=45] (n0);

      \draw[->] (n6) -- (n1);
      \draw[->] (n1) to[bend left=45] (n7);
      \draw[->] (n7) to[bend left=45] (n1);

      \draw[->] (n2) -- (n8);
      \draw[->] (n8) to[bend left=45] (n3);
      \draw[->] (n3) to[bend left=45] (n8);
    \end{tikzpicture}
  \end{minipage}%
  \hfill
  \begin{minipage}[c]{0.58\textwidth}
  \centering
  {
  \newcommand{\nodea}{\node[draw=none, inner sep=2pt] (a) {}
  ;}
  \newcommand{\nodeb}{\node[draw=none, inner sep=2pt] (b) {}
  ;}
  \newcommand{\nodec}{\node[draw=none, inner sep=2pt] (c) {}
  ;}
  \newcommand{\noded}{\node[draw=none, inner sep=2pt] (d) {}
  ;}
  \newcommand{\nodee}{\node[draw=none, inner sep=2pt] (e) {}
  ;}
  \newcommand{\nodef}{\node (f) {$\overline{11}$}
  ;}
  \newcommand{\nodeg}{\node (g) {$\overline{10}$}
  ;}
  \newcommand{\nodeh}{\node (h) {$\overline{9}$}
  ;}
  \newcommand{\nodei}{\node (i) {$\overline{8}$}
  ;}
  \newcommand{\nodej}{\node (j) {$\overline{7}$}
  ;}
  \newcommand{\nodek}{\node (k) {$\overline{6}$}
  ;}
  \newcommand{\nodel}{\node (l) {$\overline{5}$}
  ;}
  \newcommand{\nodem}{\node (m) {$\overline{4}$}
  ;}
  \newcommand{\noden}{\node (n) {$\overline{3}$}
  ;}
  \newcommand{\nodeo}{\node (o) {$\overline{2}$}
  ;}
  \newcommand{\nodep}{\node (p) {$\overline{1}$}
  ;}
  \newcommand{\nodeq}{\node (q) {$0$}
  ;}
  \newcommand{\noder}{\node (r) {$1$}
  ;}
  \newcommand{\nodes}{\node (s) {$2$}
  ;}
  \newcommand{\nodet}{\node (t) {$3$}
  ;}
  \newcommand{\nodeu}{\node (u) {$4$}
  ;}
  \newcommand{\nodev}{\node (v) {$5$}
  ;}
  \newcommand{\nodew}{\node (w) {$6$}
  ;}
  \newcommand{\nodex}{\node (x) {$7$}
  ;}
  \newcommand{\nodey}{\node (y) {$8$}
  ;}
  \newcommand{\nodez}{\node (z) {$9$}
  ;}
  \newcommand{\nodeaa}{\node (aa) {$10$}
  ;}
  \newcommand{\nodeab}{\node (ab) {$11$}
  ;}
  \newcommand{\nodeac}{\node (ac) {$12$}
  ;}
  \newcommand{\nodead}{\node (ad) {$13$}
  ;}
  \newcommand{\nodeae}{\node[draw=none, inner sep=2pt] (ae) {}
  ;}
  \newcommand{\nodeaf}{\node[draw=none, inner sep=2pt] (af) {}
  ;}
  \newcommand{\nodeag}{\node[draw=none, inner sep=2pt] (ag) {}
  ;}
  \newcommand{\nodeah}{\node[draw=none, inner sep=2pt] (ah) {}
  ;}
  \newcommand{\nodeai}{\node[draw=none, inner sep=2pt] (ai) {}
  ;}
  \begin{tikzpicture}[auto, font=\footnotesize]
  \matrix[column sep=.3cm, row sep=.5cm,ampersand replacement=\&]{
    \nodea  \&  \nodeb  \&  \nodec  \&  \noded  \&           \&  \nodee  \\
    \nodef  \&  \nodeg  \&           \&  \nodeh  \&  \nodei  \&  \nodej  \\
    \nodek  \&  \nodel  \&  \nodem  \&  \noden  \&           \&  \nodeo  \\
    \nodep  \&  \nodeq  \&           \&  \noder  \&  \nodes  \&  \nodet  \\
    \nodeu  \&  \nodev  \&  \nodew  \&  \nodex  \&           \&  \nodey  \\
    \nodez  \&  \nodeaa  \&           \&  \nodeab  \&  \nodeac  \&  \nodead  \\
    \nodeae  \&  \nodeaf  \&  \nodeag  \&  \nodeah  \&           \&  \nodeai  \\
  };

  \path[ultra thick]   (g) edge (k)
  	  (g) edge (l)
  	  (h) edge (m)
  	  (h) edge (n)
  	  (j) edge (o)
  	  (l) edge (p)
  	  (l) edge (q)
  	  (n) edge (r)
  	  (o) edge (s)
  	  (o) edge (t)
  	  (q) edge (u)
  	  (q) edge (v)
  	  (r) edge (w)
  	  (r) edge (x)
  	  (t) edge (y)
  	  (v) edge (z)
  	  (v) edge (aa)
  	  (x) edge (ab)
  	  (y) edge (ac)
  	  (y) edge (ad);
  \path[ultra thick, dotted]   (b) edge (f)
  	  (b) edge (g)
  	  (d) edge (h)
  	  (e) edge (i)
  	  (e) edge (j)
  	  (aa) edge (ae)
  	  (aa) edge (af)
  	  (ab) edge (ag)
  	  (ab) edge (ah)
  	  (ad) edge (ai);
  \end{tikzpicture}}
  \end{minipage}
        \caption{The coresolution quiver (left) and the coresolution forest (right) of the Nakayama algebra $ A $ with coKupisch series $ d = [5,4,4,5,4,5,5,6,5,4]$. The length of $ A $ is two, and its weight is $ ( 5 + 5 ) / 10 = 1 $.}
        \label{fig:coresforest-5445455654}
    \end{figure}
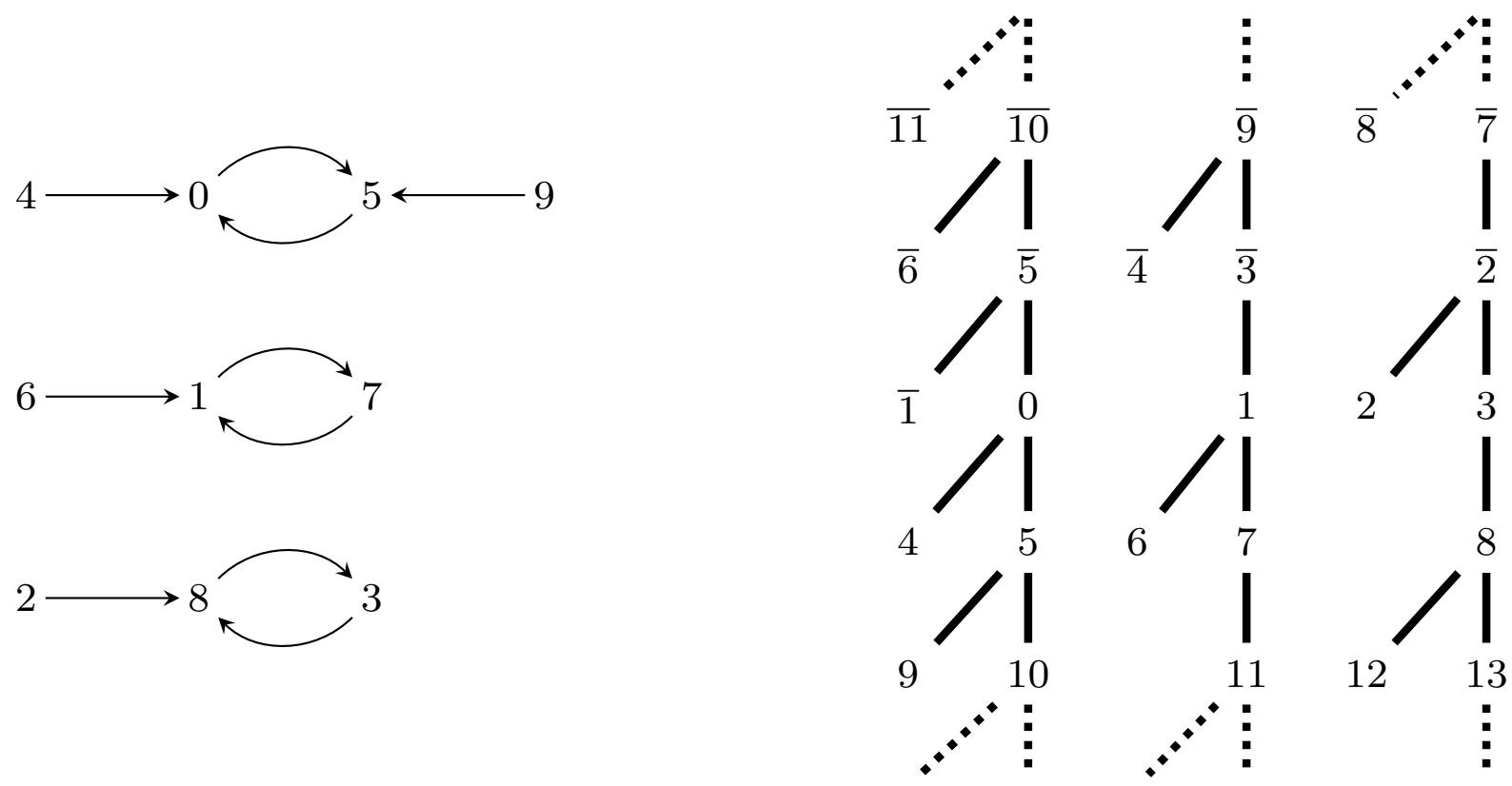
  The weight of each cycle of $ \psiquiver $ is one, and thus each component of $ \psiforest $ is a universal cover of a unique component of $ \psiquiver $ via $ \rho $.
\end{example}

We conclude this subsection by showing the equivalence of the following combinatorial and homological conditions:

\begin{proposition}\label{prop:self_injective}
  For a Nakayama algebra $ A $ and its coresolution forest $ \psiforest $, the following are equivalent:
    \begin{enumerate}
        \item\label{self-inj} $A$ is self-injective;
        \item the coKupisch series of $ A $ is constant;
        \item\label{only.vertebra}
        every vertex of $ \psiforest $ is a vertebra, i.e., $\psiforest$ has no leaves;
        \item $ \findim A = 0 $.
    \end{enumerate}
    Furthermore, if $ A $ is self-injective, then so is $ \sAus(A) $.
\end{proposition}

\begin{proof}
    $ (1) \Leftrightarrow (2) $: $ A $ is self-injective if and only if $ A $ is the field $ \K $ or cyclic and truncated (\cite[Proposition V.3.8]{ASS}).  $ A $ is the field $ \K $ if and only if $ d_i = 1 $.  $ A $ is cyclic and truncated, i.e., $ A = \K Q / I $ where $ Q $ is a cyclic quiver and $ I $ is generated by all paths of a given length $ d $, if and only if $ d_i = d $.

    $ (2) \Leftrightarrow (3) $: $\psiforest$ has no leaves if and only if $ \psi $ is surjective. Since $ \psi $ is also order preserving, this is equivalent to $\psi(i+1) - \psi(i) = 1 + d_i - d_{i+1} = 1$ for all $i \in \Z $.

    $ (1) \Leftrightarrow (4) $: Every self-injective algebra $ A $ has $ \findim A = 0 $. Now let $A$ be a non-self-injective Nakayama algebra and $P$ an indecomposable projective $A$-module that is not injective. Since $ A $ is QF-3 by \cite[Theorem~32.2]{AF}, the injective envelope $ I(P) $ of $ P $ is a projective module; see also \cite[Theorem~31.6]{AF}.
    Therefore, $ \pdim \Sigma (P) = 1 $ since the monomorphism of the short exact sequence $ 0 \to P \to I(P) \to \Sigma(P) \to 0 $ does not split, and thus $ \findim A \geq 1 $.

    The last statement is a direct consequence of \Cref{theoreminjcoresstableaus}.(\ref{theoreminjcoresstableaus3}).
\end{proof}

\subsection{Fundamental properties of coresolution forests}\label{fundamental.properties}

This subsection establishes some elementary properties of coresolution forests.
We fix an $ n $-coKupisch series $(d_i)_{i\in\Z}$ and its coresolution forest $ \psiforest $ throughout.

Let us begin with a useful new notion:

\begin{definition}
  A subset of vertices $\mathcal{I}\subseteq\psiforest$ is a \Dfn{floor} if:
  \begin{enumerate}[\rm(i)]
  \item any two vertices in $ \mathcal I $ are in the same component of $\psiforest$ and have the same distance to a common ancestor;
  \item $ \mathcal I $ is maximal with respect to that property.
  \end{enumerate}
\end{definition}

Floors satisfy the following elementary properties:

\begin{lemma}\label{prop:floors}
    The following hold for floors in $ \psiforest $:
        \begin{enumerate}
            \item Every floor is an interval.
            \item Any two floors are either equal or disjoint.
            \item The image of a floor under $ \psi $ is contained in a floor.
            \item The inverse image of a floor under $ \psi $ is a floor.
            \item\label{unique.vertebra} Every floor contains exactly one vertebra.
        \end{enumerate}
\end{lemma}

\begin{proof}
    The lemma is a direct consequence of $\psi$ being order preserving and \Cref{prop:covering}.
\end{proof}

Next, we establish a natural planar embedding of $\psiforest$, essentially unique up to the choice of a component in $\psiforest$. Most of our arguments in the rest of the paper rely on this embedding.

\begin{proposition}[Naturally labeled property]\label{thm:naturally}
  After choosing any component $T$ of the coresolution forest
  $\psiforest$, there is a uniquely determined function
  $\level:\mathbb Z\to\mathbb Z$, such that
  \begin{enumerate}
  \item\label{init} $\level(0)=0$ and the smallest vertex with level $0$ is in
    $T$,
  \item\label{parent} $\level(\psi(i)) = \level(i)-1$,
  \item\label{monotone} $\level(i+1)-\level(i)\in\{0,1\}$.
  \end{enumerate}
\end{proposition}
\begin{proof}
  Taking into account that $T$ contains a spine, \Cref{prop:floors}
  implies that there is a doubly infinite sequence
  $(\mathcal I_\ell)_{\ell\in\Z}$ of floors partitioning the vertex
  set of $T$, which we index such that the largest nonpositive
  vertex of $T$ is in $\mathcal I_0$ and
  $\psi(\mathcal I_\ell)\subseteq \mathcal I_{\ell-1}$.  Let
  $a_\ell = \min\mathcal I_\ell$.  Because $\psi$ is strictly
  decreasing, $a_\ell$ is strictly increasing.  Thus we can define
  $\level(i) = \ell$ for $i\in[a_\ell, a_{\ell+1}-1]$.

  We now check the three conditions.  \Cref{init} follows from
  $a_0\leq 0<a_1$ and the definition of $a_0$.  \Cref{monotone} is a
  consequence of $[a_\ell, a_{\ell+1}-1]$ being nonempty.  To prove
  \Cref{parent}, let $i\in [a_\ell, a_{\ell+1}-1]$.  We show that
  $a_{\ell-1}\leq\psi(a_\ell)\leq\psi(i)\leq\psi(a_{\ell+1}-1)<
  a_\ell$.

  The first inequality follows from
  $\psi(a_\ell)\in\mathcal I_{\ell-1}\subseteq [a_{\ell-1},
  a_\ell-1]$.  The next two inequalities hold because $\psi$ is order
  preserving.  We prove the last inequality by contradiction.
  Suppose that $a_\ell\leq\psi(a_{\ell+1}-1)$.  Then, since
  $\psi(a_{\ell+1}-1)\leq\psi(a_{\ell+1})$, and $a_\ell$ and
  $\psi(a_{\ell+1})$ are both in $\mathcal I_\ell$, so is
  $\psi(a_{\ell+1}-1)$.  Thus, $a_{\ell+1}-1\in\mathcal I_{\ell+1}$,
  which contradicts the choice of $a_{\ell+1}$ being the minimal
  element of $\mathcal I_{\ell+1}$.

  For uniqueness, note that the value of $\level$ on any vertex $i$
  determines $\level$ on all vertices in the same connected
  component.  Indeed, given a vertex $j$, consider any common
  ancestor $\psi^r(i) = \psi^s(j)$ of $i$ and $j$.  Then
  \Cref{parent} implies
  \[
    \level(j) = \level(\psi^s(j)) + s %
    = \level(\psi^r(i)) + s = \level(i) - r + s.
  \]
  In particular, as $ \psi ( a_{\ell+1} ) $ and $ a_\ell $ belong to the same floor $ \mathcal I_\ell $, $ \level( a_{ \ell+1 } ) = \level( a_\ell ) + 1 $ for every $ \ell \in \Z $.

    Next, we show that $ \level ( a_0 ) = 0 $. Let $ \mathcal I_0 = [ a_0 , b ] $, and recall that $ b $ is the largest nonpositive vertex in $ T $. By the previous paragraph and \Cref{init,monotone}, we have $ \level(a_0) = \level(b) \leq \level(0)=0$.
    Assuming $\level(a_0)<0$, we
  have $b<0$ by the first part of \Cref{init}.  But then the smallest
  vertex with level $0$ must be in the interval $[b+1, 0]$, which has
  no vertex in $T$, contradicting the second part of \Cref{init}.
  Therefore $ \level ( a_0 ) = 0 $, and it follows that $ \level( a_\ell ) = \ell $ for every $ \ell \in \Z $ by the previous paragraph.

  Suppose now that $\level(a_1-1)\geq1$.  Then
  $\level(\psi(a_1-1))\geq0$, and, since $a_0$ is the smallest vertex
  of level $0$ and $ \psi $ is order preserving, $a_0 \leq \psi(a_1-1) \leq \psi( a_1 )$. Since $ a_0 $ and $ \psi( a_1 ) $ are in $ \mathcal I_0 $, it follows that $ \psi( a_1-1 ) $ is also in $ \mathcal I_0 $, and thus $ a_1 - 1 \in \mathcal I_1 $, contradicting the minimality of $ a_1 $ in $ \mathcal I_1 $. Therefore, $ \level( a_1 - 1 ) = 0 $ and, by \Cref{monotone}, $ \level( j ) = 0 $ for every $ j \in [ a_0 , a_1 - 1 ] $.

    It remains to observe that every component has a vertex in $ [a_0,a_1-1] $. This follows from the facts that $ \psi $ is order preserving and every component of $ \psiforest $ contains a doubly infinite path (\Cref{prop:covering}).
\end{proof}

\begin{remark}\label{rem:component.order}
    Floors are the intersections of components with levels according to \Cref{prop:floors,thm:naturally}. More specifically, every level consists of exactly one floor from each component, and the order in which components occur repeats itself cyclically. This cyclic order does not depend on the choice of the initial component $ T $, and thus the terms \Dfn{previous/next component} are unambiguous.
\end{remark}

In what follows, the coresolution forest $\psiforest$ is assumed to be equipped with a level function as in \Cref{thm:naturally}, with respect to some arbitrary component $T$, unless stated otherwise. Furthermore, we draw $\psiforest$ so that \emph{vertices increase from left to right} on each level, and \emph{levels increase from top to bottom}. In particular, our convention that arrows are directed from bottom to top is satisfied; see the discussion preceding \Cref{ex:psiforest}.

\begin{remark}\label{levels.of.res.forests}
    The dual of \Cref{thm:naturally} holds for the resolution forest $ \gammaforest $ of an $ n $-coKupisch series $(d_i)_{i\in\Z}$. The only difference is that we require
        \begin{enumerate}
            \item[{\it(2')}]\label{antimonotone} $\level(\gamma(i)) = \level(i)+1$
        \end{enumerate}
    instead of {\it(2)} for reasons that will become apparent in \Cref{sec:main-results-symmetry}; see \Cref{minus.forest}. Pictorially, we draw $\gammaforest$ so that \emph{vertices increase from right to left} on each level, and \emph{levels increase from bottom to top}, ensuring that arrows are directed from bottom to top as in the coresolution forest.
\end{remark}

The next lemma, which locates $\gamma(i) - 1 $ for any vertex $i \in \psiforest $, will play a significant role throughout the paper.
Note that $ \gamma(i)-1 $ has a natural homological interpretation as the top of the indecomposable injective $ A $-module corresponding to $ i $ when $ \psiforest $ is the coresolution forest of a Nakayama algebra $ A $.

\begin{lemma}[$S$-Lemma]
    \label{lemma:largest:child:non:leaf}
    For every vertex $ i \in \psiforest $, $\gamma(i)-1$ is the largest child of the
    largest non-leaf smaller than or equal to $i-1$.
\end{lemma}
\begin{proof}
  We first show that the parent of $\gamma(i)-1$ is smaller than or equal
  to $i-1$, and that the parent of $\gamma(i)$ is strictly larger
  than $i-1$. By \Cref{lem:kupisch-cokupisch}, the Kupisch series and the coKupisch series of $A$ are related by
  \[
    c_i = \max\{k \in \Z_{\geq1} \mid k\leq d_{i+k-1}\}.
  \]
  In particular, we have $c_i\leq d_{i+c_i-1}$, and therefore
  \[
    \psi(\gamma(i)-1)=\psi(i+c_i-1) = i+c_i-1-d_{i+c_i-1}\leq i-1.
  \]
  Similarly, $c_i+1 > d_{i+c_i}$, and therefore
  \[
    \psi(\gamma(i))=\psi(i+c_i) = i+c_i-d_{i+c_i} > i - 1.
  \]

  Since $\gamma(i)-1$ and $\gamma(i)$ are consecutive and $\psi$ is order preserving, any vertices strictly between $\psi(\gamma(i)-1)$ and $\psi(\gamma(i))$ are leaves, and $\gamma(i)-1$ is the largest child of its parent. Since $\psi(\gamma(i)-1) \leq i - 1$ and $\psi(\gamma(i)) > i - 1$, the largest vertex smaller than or equal to $i-1$ that has a child is $\psi(\gamma(i)-1)$.
\end{proof}

The following notion, which corresponds to the natural notions of left and right in ordered trees, underlies most of our main results and their proofs:
\begin{definition}\label{to.the.left}
  A vertex $i \in \psiforest $ is \Dfn{to the left} of a vertex $j$
  if $\psi^\ell(i) < \psi^m(j)$ for some pair of nonnegative integers
  $\ell$ and $m$ such that $\psi^\ell(i)$ and $\psi^m(j)$ are on the
  same level.  A vertex $j$ is \Dfn{to the right} of a vertex $i$ if
  $i$ is to the left of $j$. Similarly, a subset of vertices $ \mathcal S \subseteq \psiforest $ is to the left/right of a subset $ \mathcal T \subseteq \psiforest $ if every vertex in $ \mathcal S $ is to the left/right of every vertex in $ \mathcal T $.
\end{definition}

The above relation is transitive, an important property that can be easily verified:

\begin{lemma}
    For three vertices $i,j,k \in \psiforest $, if $i$ is to the left of $j$, and $j$ is to the left of $k$, then $i$ is to the left of $k$.
\end{lemma}

Next, we collect several useful consequences of \Cref{prop:floors,thm:naturally,lemma:largest:child:non:leaf}, which will recur in the proofs of our main theorems. For that, we need one more notion:

\begin{definition}\label{defn:left.right}
  A vertex of $ \psiforest $ is a \Dfn{right vertex} if it lies to the right of the spine
  in its component, and a \Dfn{left vertex} is defined analogously.
\end{definition}

\begin{corollary}\label{S-lemma:corollaries}
    Let $ \psiforest $ be the coresolution forest of an $ n $-coKupisch series $(d_i)_{i\in\Z}$ and $(c_i)_{i\in\Z}$ the associated $ n $-Kupisch series. Then, for every vertex $ i \in \psiforest $:
    \begin{enumerate}
    \item\label{S-lemma:no:ancestors} No two distinct vertices in the interval $[i,\gamma(i)-1]$ are related by ancestorship in $ \psiforest $.
      \item\label{S-lemma:levels}
      For every positive integer $ k \leq c_i $, the vertices $i-1$, $i+k-1$ and $\gamma(i)-1$ are within two levels of $ \psiforest $.
     \item\label{lemma:right.vertex} If $ i-1 $ is a right vertex or a vertebra in $ \psiforest $, then $ \gamma(i)-1 $ is one level below $ i-1 $ and in the same component.
    \item\label{S-lemma:interval}
    If $i-1$ is a leaf of $ \psiforest $, then
      $[i, \gamma(i)-1]$
      consists of the vertices on level
      $\ell=\level(i-1)$ that are to the right of $i-1$, together
      with the vertices on level $\ell+1$ to
      the left of $i-1$.
    \item\label{S-lemma:ancestors} Any vertex larger than $i-1$ has a unique ancestor in the interval $[i,\gamma(i)-1] $ in $ \psiforest $.
    \end{enumerate}
\end{corollary}

\begin{proof}
    (1) If $j\leq\gamma(i)-1$, then $ \psi(j)\leq\psi(\gamma(i)-1)\leq i-1 $ due to \Cref{lemma:largest:child:non:leaf} and $ \psi $ being order preserving.
    In particular, $\psi(j)$ and $j$ cannot be both in $[i, \gamma(i)-1]$.

    (2) If $ \gamma(i)-1 $ is on the same level as $ i-1 $, then all three vertices have the same level $ \ell = \level(i-1) $ by \Cref{monotone} of \Cref{thm:naturally}. Otherwise, $ \level( \gamma(i)-1 ) > \ell $ by the same property. If $ \level ( \gamma(i)-1 ) > \ell+1 $, then \Cref{parent} of \Cref{thm:naturally} implies that $ \level ( \psi( \gamma(i)-1 ) ) > \ell $. But then $ i-1 < \psi( \gamma(i)-1 ) $ contradicting \Cref{lemma:largest:child:non:leaf}.

    (3) Since $ i-1 $ is a right vertex or a vertebra, there is a vertebra $ j \leq i-1 $ to its left and on the same level. It follows that $ \psi ( \gamma(i)-1 ) \in [j,i] $ by \Cref{lemma:largest:child:non:leaf}, implying that $ \psi ( \gamma(i)-1 ) $ is on the same floor as $ i-1 $.

    (4) If $ \gamma(i)-1 $ is on the same level as $ i-1 $, then it follows from \Cref{lemma:largest:child:non:leaf} that $ \gamma(i)-1 $ is the largest vertex of level $ \ell $. Furthermore, every vertex on level $ \ell $ to the left of $ i-1 $ is a leaf and, thus, no vertex on level $ \ell+1 $ (or greater) is to the left of $ i-1 $.
    Due to \Cref{S-lemma:levels}, the remaining possibility is that $ \gamma(i)-1 $ is one level below $ i-1 $. Then the vertices of level $ \ell+1 $ that are at most equal to $ \gamma(i)-1 $ are exactly the vertices of that level that are to the left of $ i-1 $.

    (5) If $i - 1$ and $\gamma(i) - 1$ lie on the same level of $\psiforest$, then the following two facts follow from \Cref{lemma:largest:child:non:leaf}:
    Every vertex in $[ \psi(\gamma(i)-1)+1, i-1 ]$ is a leaf, and $\gamma(i)-1$ is the largest vertex on level $\ell = \level (i-1)$. Therefore, every vertex on level $\ell$ is either in $[ i-1,\gamma(i)-1 ]$ or a leaf.
    If $ \level(\gamma(i)-1) = \ell+1 $, then $ j $ or $\psi(j) $ is in $ [ i-1, \gamma(i)-1 ] $ for every vertex $ j $ with $ \level(j) = \ell+1 $, since $ j \geq \gamma(i) $ implies $ \psi(j) \geq i $ by \Cref{lemma:largest:child:non:leaf}.
\end{proof}

We close this subsection by introducing the notion of distance and showing that it is monotone on each argument for vertices on the same level of the coresolution forest.

\begin{definition}
    The \Dfn{distance} between two vertices $i , j \in \psiforest$, denoted by $\dist(i,j)$, is the length of the unique undirected path between them if such a path exists, and $\infty$ otherwise.
\end{definition}

The distance and the \Dfn{closest common ancestor} of two vertices $ i , j \in \psiforest $, denoted by $ i \vee j $,  are connected via the level function as follows:

\begin{lemma}\label{prop:distance}
    For two vertices $ i , j $ in the same component of $ \psiforest $,
    \[
        \dist(i, j) = \level(i) + \level(j) - 2 \level(i\vee j).
    \]
\end{lemma}

\begin{proof}
    Direct consequence of the definitions and the properties of the level function as in \Cref{thm:naturally}.
\end{proof}

\begin{lemma}\label{lemma:common-ancestor}
  Let $x$, $y$ and $z$ be vertices of $\psiforest$ such that $x$ is
  left of $y$, which in turn is left of $z$.
    \begin{enumerate}
        \item If $x$, $y$ and $z$ are in the same component of
  $\psiforest$, then
  \[
    \level(x\vee z) = \min\{\level(x\vee y), \level(y\vee z)\}.
  \]
       \item If $x$ and $y$ are on the same level, then
  $\dist(x, z)\geq\dist(y, z)$.
  Similarly, if $y$
  and $z$ are on the same level, then $\dist(x, z)\geq\dist(x, y)$.
    \end{enumerate}
\end{lemma}
\begin{proof}
  (1) Assume without loss of
  generality that $\level(x\vee y)\leq\level(y\vee z)$.  Since
  $x\vee y$ and $y\vee z$ are both ancestors of $y$, $x\vee y$ is an
  ancestor of $y\vee z$ and, therefore, also an ancestor of $z$.  In
  particular, $x\vee y$ is an ancestor of $x\vee z$.  Conversely,
  since $x$ is left of $y$ and $y$ is left of $z$, a common ancestor
  of $x$ and $z$ is also an ancestor of $y$.  Thus $x\vee z$ is an
  ancestor of $x\vee y$.  We deduce that $x\vee y = x\vee z$.

  (2) If $x$, $y$ and $z$ are in the same component of $\psiforest$, and $x$ and $y$ are on the same level, then
  \Cref{prop:distance} yields $\dist(x, z)\geq\dist(y, z)$.
  If $y$ and $z$ are in distinct components, also $x$ and
  $z$ must be in distinct components, because $x$ is to the
  left of $y$.  Finally, if $x$ and $z$ are in distinct connected
  components, we trivially have $\dist(x, z)\geq\dist(y, z)$.

  The case where $y$
  and $z$ are on the same level is analogous and thus omitted.
\end{proof}

\subsection{Three-distance formula}\label{subsec:three-distance}

The goal of this subsection is to give a combinatorial way of computing the injective dimensions of projective $ \sAus(A) $-modules, where $ A $ is a Nakayama algebra. As a first application, we show that each such dimension is either even or infinite.

We begin by expressing the injective dimension of an indecomposable $ A $-module in terms of the distance between two vertices in its coresolution forest:

\begin{lemma}\label{lemma:distances}
  Let $M = \mik$ be an indecomposable module over a Nakayama algebra $ A $. Then
  \[
    \idim M = \dist(i - 1, i + k - 1) - 1,
  \]
  where the distance is computed in the coresolution forest of $ A $.
\end{lemma}

\begin{proof}
    Consider a short exact sequence
        \[
            0 \to M \to I^0 \to \Sigma (M) \to 0
        \]
    where $ M \to I^0 $ is an injective envelope of $ M = \mik $. Then, one calculates that
        \[
            I^0 = M ( \psi(i+k-1)+1 , i+k-1 - \psi(i+k-1) )
        \]
    and thus
        \[
            \Sigma (M) = M ( \psi(i+k-1)+1 , i-1 - \psi(i+k-1) ) .
        \]

    Now consider the minimal injective coresolution of $M$
    \[
    0 \rightarrow M \rightarrow I^0 \rightarrow I^1 \rightarrow I^2 \rightarrow \cdots  .
    \]
    An inductive application of the above calculations gives
        \[
            \Sigma^{ 2m } (M) = M( \psi^m(i-1)+1 , \psi^m(i+k-1)-\psi^m(i-1) )
        \]
    and
        \[
            \Sigma^{ 2m+1 } (M) = M( \psi^{m+1}(i+k-1)+1 , \psi^m(i-1)-\psi^{m+1}(i+k-1) )
        \]
    for every integer $ m \geq 0 $. As $ \idim M = d < \infty $ if and only if $ d \geq 0 $ is minimal such that $ \Sigma^{d+1} (M) =0$, we deduce that
        \[
            \idim M =
                \begin{cases}
                  2m \textrm{ (even)}, & \text{if } \psi^m(i-1) = \psi^{m+1}(i+k-1) \\
                  2m+1 \textrm{ (odd)},  & \text{if } \psi^{m+1}(i-1) = \psi^{m+1}(i+k-1) \\
                  \infty,  & \text{if } i-1, i+k-1 \textrm{ are in different components of $ \psiforest$},
                \end{cases}
        \]
    where $ m \geq 0 $ is minimal with the specified property. This confirms the desired equality in every case, as $ \dist(i-1,i+k-1) = 2m+1 $, $ 2m+2 $ or $ \infty $, respectively.
\end{proof}

\Cref{lemma:distances} yields the following characterization of Nakayama algebras with infinite global dimension, which recovers \cite[Proposition 1.1]{S2} via \Cref{rem:weight}:
\begin{proposition}\label{lemma:infinitegldim}
    Let $A$ be a Nakayama algebra.
    Then $\gldim A < \infty$ if and only if the coresolution forest of $ A $ is connected.
\end{proposition}
\begin{proof}
  If $ n $ is the number of vertices in the quiver of $ A $, then
        \[
            \gldim A = \sup \bigl\{ \dist(i-1,i) \mid i = 0 ,1 ,\ldots , n-1 \bigr\}
        \]
  since the global dimension is attained as the injective dimension of a simple module and in view of \Cref{lemma:distances}.
  Now, $\dist(i - 1, i)$ is always finite if the coresolution forest is connected.
  Conversely, if the coresolution forest is not connected, then we can pick $i$ so that $i-1$ is in a different component, and thus their distance is infinite.
  By the $ n $-periodicity of the coresolution forest, we may choose $ i $ to be in the interval $ [ 0 , n-1 ] $.
\end{proof}

The following consequence of \Cref{lemma:infinitegldim} provides a useful criterion for infinite global dimension in Nakayama algebras:

\begin{corollary}\label{cor:finite.gldim}
   Let $ A $ be a Nakayama algebra and $ \psiforest $ its coresolution forest. If $ i-1 $ and $ \gamma(i)-1 $ lie in the same floor of $ \psiforest $ for some $ i $, then $ \gldim A < \infty $.
\end{corollary}

\begin{proof}
    Since $ i-1 $ and $ \gamma(i)-1 $ lie in the same floor of $ \psiforest $, every vertex in $ [ \psi (\gamma(i)-1)+1 , i-1 ] $ must be a leaf according to \Cref{lemma:largest:child:non:leaf}. Since $ \psi(\gamma(i)-1) $ is one level above $ i-1 $ and every floor contains a (unique) vertebra, the level of $ i-1 $ comprises a unique floor. Therefore, $ \psiforest $ is connected (see \Cref{rem:component.order}) and the result follows by \Cref{lemma:infinitegldim}.
\end{proof}

We prove next a central result of this paper, calculating the injective dimensions of projective modules over the stable Auslander algebra of a Nakayama algebra $ A $ in terms of distances in its coresolution forest:

\begin{proposition}[Three-distance formula]
    \label{three:distance:formula}
    Let $A$ be a Nakayama algebra and $M $ an indecomposable non-projective $A$-module. If $ d = \min \{ \idim \Omega(M) , \idim P(M) ,\allowbreak \idim M \} $, then
    \[
      \stabidim=
      \begin{cases}
        3d-1, & \textrm{if } d = \idim \Omega(M),\\
        3d,   & \textrm{if } d = \idim P(M) < \idim \Omega(M) ,\\
        3d+1, & \text{otherwise},
      \end{cases}
    \]
    where $ \pstab M $ is the projective indecomposable $ \sAus(A) $-module corresponding to $ M $.
    
    Moreover, the injective dimensions can be computed via distances in the coresolution forest of $ A $ as follows
  \begin{align*}
      \idim \Omega(M) &= \dist(i+k-1, \gamma(i)-1)-1, \\
      \idim P(M)      &= \dist(i-1, \gamma(i)-1)-1, \\
      \idim M         &= \dist(i-1, i+k-1)-1,
   \end{align*}
    where $ M = \mik $ has top isomorphic to $ S_i $ and Loewy length $ k < c_i $.
\end{proposition}

\begin{proof}
    Over a Nakayama algebra $A$, indecomposable modules are closed under projective covers, syzygies, injective envelopes and cosyzygies. As a first consequence, we  can compute the injective dimensions of $ M $, $\Omega(M) = \mik[i + k][c_i - k]$ and $P(M) = \mik[i][c_i]$ using \Cref{lemma:distances}.

       We observe next that the injective coresolution of $ \pstab M $ in \Cref{theoreminjcoresstableaus} (\ref{theoreminjcoresstableaus3}) is minimal after truncation at the first zero term. Indeed, all its nonzero injective $ \sAus(A) $-modules are indecomposable by dimension shifting and \Cref{injresremark}.(\ref{rem:injective.sAus}), and the maps between consecutive nonzero Ext spaces before the first zero term are nonzero. The latter follows from uniseriality and dimension shifting, using the criterion that $\Ext_A^1(-,f)=0$ if and only if $f$ factors through an injective module; see \cite[discussion following Proposition~I.3.3, p.~18]{AR}. For the construction of this coresolution from the long exact Ext-sequence, see \cite[proof of Proposition~I.3.2, pp.~16--17]{AR}. The result follows now from the fact that, for every $ A $-module $ N $, $\Ext^k(X, N) = 0$ if and only if $\idim N \leq k-1$.
\end{proof}

Observe that, when $ \idim \pstab M$ is finite, its residue modulo $ 3 $ reveals information about how $ \idim \Omega(M)$, $ \idim P(M) $ and $ \idim M $ relate to each other. In particular, it reveals one of the modules where the minimum injective dimension is attained.

\begin{example}\label{ex:3.dist.form}
  We revisit the Nakayama algebra $ A $ with coKupisch series $d = [1,2,3,4,4,5,4]$ of \Cref{example:linearnakayama}; see \Cref{fig:coresforest-1234454} for its coresolution forest.
  Three triplets of vertices corresponding to the modules $ \mik[2][1] , \mik[3][1] $ and $ \mik[2][3] $ are highlighted in different colors.
  We have also highlighted the paths between the pairs in each triplet that give the minimum distance for the three-distance formula.
  Observe that, in the three-distance formula for a module $ M = \mik $, we may use any $ i' $ instead of $ i $ as long as they are equivalent modulo the number $ n $ of all vertices (here $ n = 7 $).
  We collect the relevant data in the following table:
  \smallskip
  \begin{center}
    {\small\setlength{\tabcolsep}{4pt}
    \begin{tabular}{c|ccc|ccc|c|c}
      $M$ & $i-1$ & $i+k-1$ & $\gamma(i)-1$ & $\idim \Omega(M)$ & $\idim P(M)$ & $\idim M$ & $d$ & $\idim \pstab{M}$ \\
      \hline
      $\mik[2][1]$ & $-6$ & $-5$ & $-2$ & 2 & 2 & \textbf{1} & 1 & 4 \\
      $\mik[3][1]$ & $2$  & $3$  & $6$  & 2 & \textbf{0} & 1 & 0 & 0 \\
      $\mik[2][3]$ & $8$  & $11$ & $12$ & \textbf{1} & 2 & 2 & 1 & 2 \\
    \end{tabular}}
  \end{center}
\end{example}

To compute the minimal distance $ d+1 $ in \Cref{three:distance:formula}, we climb up the tree from each of the three vertices $\gamma(i)-1$, $i+k-1$ and $i-1$ associated to $ M = \mik $ in turn, until two of the resulting paths meet:
\begin{definition}
    For an indecomposable non-projective module $M = \mik$ over a Nakayama algebra $A$, the \Dfn{climbing sequence} $\{a_j\}_{j \geq 0}$ of $M$ is defined using the rules $a_0 = \gamma(i) - 1$, $a_1 = i + k - 1$, $a_2 = i - 1$, and $a_j = \psi(a_{j - 3})$ for $j \geq 3$.
\end{definition}

\begin{proposition}\label{prop:even}
  Let $M = \mik$ be an indecomposable non-projective module over a Nakayama algebra $ A $. Then
    \[\stabidim = 2(\ell - 2),\]
  where $\ell \geq 0$ is the minimum index such that $a_\ell = a_{\ell + 1}$ in the climbing sequence. If no such $\ell$ exists, then $\stabidim = \infty$.
  In particular, $\stabidim$ is either even or infinite.
\end{proposition}
\begin{proof}
  Let $d+1$ be the minimal distance between $i-1$, $i+k-1$ and
  $\gamma(i)-1$, as in \Cref{three:distance:formula}.
  Assume first
  that $d = \idim\Omega(M) = \dist(i+k-1, \gamma(i)-1) - 1$.
  We
  show that in this case $i+k-1$ and $\gamma(i)-1$
  must be on the
  same level.
  Suppose the contrary.
  Then $\psi(\gamma(i)-1)$ is
  left of $i-1$, which in turn is left of $i+k-1$ and all three are
  on the same level.
  In this situation \Cref{lemma:common-ancestor} implies
  $\dist(i-1, i+k-1)\leq\dist(\psi(\gamma(i)-1), i+k-1) <
  \dist(\gamma(i)-1, i+k-1)$, contradicting the assumption
  $d=\idim\Omega(M)$.

  When $i+k-1$ and $\gamma(i)-1$ are on the same level, their distance $d+1$ is even.  More precisely, let $m$ be minimal such that
  $\psi^m(\gamma(i)-1) = \psi^m(i+k-1)$. Then $d = 2m-1$ and $\ell = 3m$.
  \Cref{three:distance:formula} implies $\stabidim = 3d-1$, so $\stabidim = 2(\ell-2)$.

  Assume now that
  $d = \idim P(M) = \dist(i-1, \gamma(i)-1) - 1 < \idim\Omega(M)$.
  Similar to the reasoning above, we find that $i-1$ and
  $\gamma(i)-1$ must be on different levels.
  Otherwise, all three of
  $i-1$, $i+k-1$ and $\gamma(i)-1$ are on the same level, and
  \Cref{lemma:common-ancestor} implies
  $\dist(i+k-1, \gamma(i)-1)\leq\dist(i-1, \gamma(i)-1)$, a
  contradiction.

  In this case, let $m$ be minimal such that
  $\psi^m(i-1) = \psi^{m+1}(\gamma(i)-1)$.
  Then $d = 2m$ and $\ell=3m+2$, so $\stabidim = 3d = 2(\ell-2)$.

  Finally, assume that
  $d = \idim M = \dist(i-1, i+k-1) - 1 < \idim P(M)$ and $\idim M < \idim\Omega(M)$.
  We find that
  $i-1$ and $i+k-1$ must be on the same level, because otherwise
  \Cref{lemma:common-ancestor} implies that
  $\dist(\gamma(i)-1, i-1)\leq\dist(i-1, i+k-1)$, contradicting
  $\idim M < \idim P(M)$.

  Let $m$ be minimal such that $\psi^m(i-1) = \psi^m(i+k-1)$, so $d=2m-1$. Then $\ell=3m+1$ and \Cref{three:distance:formula} implies $\stabidim = 3d + 1 = 2(\ell-2)$.
\end{proof}

\begin{example}
     Let $ A $ be the Nakayama algebra of \Cref{ex:fig.2} with coKupisch series $ d = [8,9,9,9] $ and coresolution forest as in \Cref{fig:coresforest-8999}. The climbing sequence of the module $M=\mik[0][4]$ is $$19,15,11,10,6,2,1,-3,-7,-8,-12,-16,-16,\dots $$
     implying that $\ell = 11$ and $\idim\pstab{M} = 2(11-2) = 18$.
\end{example}

\subsection{Height and the finitistic dimension}\label{compute.fin.dim}

The aim of this subsection is to connect the height of a Nakayama algebra, a notion that goes back to \cite{G}, with its finitistic dimension. In particular, we calculate the exact value of the finitistic dimension of every Nakayama algebra.

\begin{definition}\label{def:height}
    The \Dfn{height of a vertex $i $ of $ \psiforest $}, denoted by $\height(i)$, is the distance of $ i $ to a spine, that is, the number of edges in the shortest path from $ i $ to a vertebra. 
    The \Dfn{height of $\psiforest$}, denoted by $\height(\psiforest)$, is the maximum height of any vertex in the forest. If $ \psiforest $ is the coresolution forest of a Nakayama algebra $ A $, then we call this the \Dfn{height of $ A $} and denote it by $ \height(A) $. Furthermore, a vertex $ i $ such that $ \height(i) $ attains this maximum will be called a \Dfn{vertex of maximum height}.
\end{definition}

\begin{remark}
    Every vertex of maximum height is necessarily a leaf since, otherwise, the height of its children would violate this maximality.
\end{remark}

The height of a Nakayama algebra $ A $ is not a new notion. It goes back at least to \cite{G}, since $ \height(A) $ is equal to the maximum distance of a vertex in $\psiquiver$ to a cyclic vertex by \Cref{prop:covering}.

We also note that the height is a left-right symmetric notion by \cite[Corollary~2, Appendix~A.1]{Ri2}, i.e., $ \height (A) = \height (A^{\op}) $ for every Nakayama algebra $ A $.

We begin by establishing a first connection between the height and Iwanaga--Gorensteinness:\begin{lemma}\label{lem:non.IG}Let $A$ be a Nakayama algebra and $\psiforest$ its coresolution forest. If $\gldim A = \infty$ and there exists a vertex $i-1 \in \psiforest$ of maximum height that is also left, then $A$ is not Iwanaga--Gorenstein. In particular, $\idim P_i = \dist(i-1,\gamma(i)-1) = \infty$.\end{lemma}

\begin{proof}
    If $i-1$ and $ \gamma(i)-1 $ are on the same level, then their distance is infinite by \Cref{cor:finite.gldim}. Otherwise, $ \level(\gamma(i)-1) = \level(i-1)+1 $ by \Cref{S-lemma:corollaries}.(\ref{S-lemma:levels}).
    Since $ i-1 $ is a left leaf of height $ \height(A) $, every vertex $ j $ to its left and on the same floor must also be a leaf (of height $ \height(A) $), otherwise any child of $j$ would violate maximality of $ \height(i-1) $. Therefore, $ \gamma(i)-1 $ is in a different component than $i-1$ in view of \Cref{lemma:largest:child:non:leaf}.
    It remains to recall that $ P_i = M( i, c_i ) $, and thus $ \idim P_i = \dist ( i-1, i+c_i-1 ) = \dist( i-1, \gamma(i)-1 ) = \infty $ by \Cref{lemma:distances} and the definition of $ \gamma $.
\end{proof}

In order to establish a homological interpretation of the height, we need the following notion, which will also play a role in the next section:

\begin{definition}
     A vertex is \Dfn{rightmost} (respectively, \Dfn{leftmost}) with respect to some property if there is no vertex to its right (respectively, left) and on its floor with this property. When no property is specified, then these terms are meant with respect to all vertices on the same floor.
\end{definition}

Note that rightmost and leftmost vertices may be on a spine, unlike right and left vertices introduced in \Cref{defn:left.right}.

\begin{proposition}\label{prop:findimA}
    For a Nakayama algebra $ A $, we have
         \[
            \findim A =
            \begin{cases}
                2 \height(A) , & \text{if $ \psiforest $ has a non-left vertex of maximum height} \\
                2 \height(A)-1  , & \text{otherwise.}
                    \end{cases}
        \]
\end{proposition}

\begin{proof}
    We are going to use the equalities
        \[
            \findim A = \findim A^\op = \sup \bigl\{ \idim M \mid M\in\mod A \text{ with }\idim M<\infty \bigr\}    .
        \]
    The former is due to \cite[Theorem~A.(vii)]{Se2} for cyclic Nakayama algebras and to the left-right symmetry of the global dimension for linear Nakayama algebras, while the latter is due to the standard duality between left and right finitely generated $ A $-modules.
    
    Let $M = \mik$ be an indecomposable $ A $-module such that $ \idim M < \infty $, where $ \idim M = \dist(i-1, i+k-1) - 1 $ according to \Cref{lemma:distances}. Since the vertices $ i-1 $ and $ i+k-1 $ are within two consecutive levels of $ \psiforest $ (\Cref{S-lemma:corollaries}.(\ref{S-lemma:levels})), we get $ \dist(i-1, i+k-1) \leq 2 \height(A) + 1 $ due to \Cref{cor:distance}, and thus $ \idim M \leq 2 \height(A) $. We conclude that  $ \findim A^{\op} \leq 2 \height ( A ) $.

    Next, we show that $ \findim A^{\op} \geq 2 \height ( A ) - 1 $. We assume that $ \height ( A ) > 0 $, as otherwise the inequality is trivial, and construct a simple $ A $-module $ S_i $ with $ \idim S_i = 2 \height ( A ) - 1 $. Let $ j $ be a leftmost vertex that is right and of maximum height. Since the height of $ j-1 $ is less than the height of $ j $, the vertex $ j-1 \vee j $ is a vertebra and the path connecting $ j-1 $ to $ j-1 \vee j $ contains at least one more vertebra. Therefore, $ \dist (  j-1, j ) = 2 \height ( A ) $ and we set $ i = j $. If there is no vertex that is right and of maximum height, there has to be a left vertex $ j $ that is of maximum height. Setting $ i=j+1 $ for a rightmost such $ j $ and arguing as before gives the claimed simple module.

    Now assume that $ \findim A^{\op} = 2 \height (A) $, and take $ M = \mik $ such that $ \idim M = 2 \height (A) $. Then $ \dist(i-1, i+k-1) = 2 \height(A) + 1 $, implying that $ i+k-1 $ is one level below $ i-1 $, and $ i-1 $ is $ \height(A) $ levels below $ i-1 \vee i+k-1 $. In particular, the distance between $ i+k-1 $ and $ i-1 \vee i+k-1 $ is $ \height(A)+1 $, implying that the last vertex before $ i-1 \vee i+k-1 $ in the path connecting $i+k-1 $ to $ i-1 \vee i+k-1 $ is a vertebra. Therefore, $ i-1 $ is a vertex that is right and of maximum height.

    Conversely, let $ i-1 $ be a vertex that is right and of maximum height, and we may assume that $ i-1 $ is a rightmost such vertex. We complete the proof by showing that $ \idim P_i = 2 \height(A) $ forcing $ \findim A^{\op} = 2 \height(A) $. Recall that $ P_i = M (i,c_i) $ and $ i+c_i-1 = \gamma(i)-1 $. The vertex $ i-1 \vee \gamma(i)-1 $ is a vertebra as the vertex $ \gamma(i)-1 $ is one level below $ i-1 $ due to \Cref{S-lemma:corollaries}.(\ref{lemma:right.vertex}); otherwise, we would have $ \height ( \gamma(i)-1 ) = \height ( i-1 ) + 1 $, a contradiction to the maximality of $ \height ( i-1 ) $. We deduce that $ \dist( i-1, i+c_i-1 ) = 2 \height(A) + 1 $, implying that $ \idim P_i =  \dist( i-1, i+c_i-1  ) - 1 = 2 \height(A) $ as wanted.
\end{proof}

\begin{remark}
    In the proof of \Cref{prop:findimA}, we used the fact that $ \findim A = \findim A^\op $ for every Nakayama algebra $ A $, a fact due to Sen \cite{Se2} which we recover in \Cref{cor:findim.findimop}.
\end{remark}

\Cref{prop:findimA} yields the following homological interpretation of the height:

\begin{corollary} \label{prop:heightfindim}
For every Nakayama algebra $A$,
\[
  \height(A)=\ceilfindim{A}.
\]
\end{corollary}

\section{Proofs of main results}\label{sec:main-results}

This section is devoted to the proofs of the main results.

\subsection{The finitistic Gorenstein dimension}\label{sec:main-results-bounds}

The main goal of this subsection is to prove \Cref{mainresult1} from the introduction. Specifically, we compute the finitistic Gorenstein dimension of $\sAus(A)$ in terms of the height of the Nakayama algebra $ A $, provided that the global
dimension of $A$ is infinite; see \Cref{prop:infinite:height}.
We then translate this computation via \Cref{prop:findimA} to get an expression in terms of the finitistic dimension of $ A $; see \Cref{cor:finGordim&findim}. We also provide a direct proof of \Cref{lettercor}; see \Cref{prop:iwanaga-gorenstein}.

Finally, we show that the injective dimension of the projective $ \sAus(A) $-module $ \pstab S $ corresponding to any simple non-projective $ A $-module $ S $ is always finite, regardless of whether $ \gldim A $ is finite; see \Cref{thm:idim-simple-finite}. In conjunction with \Cref{prop:infinite:height}, this yields part of \Cref{mainresult2} (namely, the case where $\gldim A = \infty$).

We begin with a useful preliminary result:

\begin{lemma}\label{cor:distance}
    For two vertices $ i $ and $ j $ in the same component of $ \psiforest $,
        \[
            \dist(i,j) \leq 2 \height(A) + | \level(i) - \level(j) | .
        \]
    Furthermore, if $ \level(i) < \level(j) $, equality is attained if and only if the height of $ i $ is maximum.
\end{lemma}

\begin{proof}
    If $ \height(A) = 0 $, then the result follows from \Cref{prop:self_injective}. We may therefore assume that $ \height(A) > 0 $.

    Since every floor contains exactly one vertebra (\Cref{prop:floors}), at least one of the paths leading from $ i$ or $ j $ to $ i \vee j $ does not contain a vertebra, except possibly for $ i \vee j $ itself. Let $ x $ be such a vertex, and $ y $ the remaining one.
    Then $ d = \dist(x,i\vee j) \leq \height(x) \leq \height(A) $, and $ \dist(y,i\vee j) = d \pm | \level(i) - \level(j) | $ depending on whether $ y $ is above $ x $ or not. All in all, $ \dist(i,j) = 2d \pm | \level(i) - \level(j) | \leq 2 \height(A) + | \level(i) - \level(j) | $.

    If $ \dist(i,j) = 2 \height(A) + | \level(i) - \level(j) | $, then $ d = \height (A) $ and $ \level(x) \leq \level(y) $. Therefore, $ \height (A) = d \leq \height(x) \leq \height(A) $ implies that $ x $ is of maximum height.
    Furthermore, if $ \level(i) < \level(j) $, then $ x=i $, proving one implication.
    Conversely, assume that $ \level(i) < \level(j) $ and $ i $ is of maximum height. Then $ i\vee j $ is a vertebra as, otherwise, $ \height(j) $ would violate the maximality of $ \height(i) $. Furthermore, $i\vee j $ has to be the closest vertebra to $ i $, implying that $\dist(i,j) = 2 \height(A) + | \level(i) - \level(j) | $.
\end{proof}

Next, we establish an upper bound for the injective dimensions of certain projective $ \sAus(A) $-modules:
\begin{lemma}\label{lemma:idim:bound:same:tree}
    Let $A$ be a Nakayama algebra and $\mik$ an indecomposable non-projective $A$-module such that $i-1$, $i+k-1$, and $\gamma(i) - 1$ all lie in the same component of $\psiforest$.
    Then \[\stabidim \leq 6\height(A) - 2.\]
\end{lemma}
\begin{proof}
    Let $d+1$ be the minimal distance between any two of the three vertices $i-1$, $i+k-1$, and $\gamma(i) - 1$.  Then \Cref{three:distance:formula} guarantees that $\stabidim\leq 3d+1$.

    By \Cref{S-lemma:corollaries}.(\ref{S-lemma:levels}), at least two vertices among $i-1$, $i+k-1$, and $\gamma(i) - 1$ must be on the same level.  Therefore, their distance is at most $2\height(A)$ due to \Cref{cor:distance}.  It follows that $\stabidim \leq 3d+1\leq 3(2\height(A) - 1) + 1 = 6\height(A) - 2$.
\end{proof}

Digressing shortly from the main goal of this subsection, we show that \Cref{lemma:idim:bound:same:tree} already yields an optimal bound of $ \gldim \sAus(A) $ for Nakayama algebras $ A $ with $ \gldim A < \infty $:
\begin{proposition}\label{prop:finite:height}
    For Nakayama algebras $A$ with $\gldim A < \infty$,
    \[
    \gldim\sAus(A)
    \leq 6\height(A) - 2,
    \]
    and the bound is optimal.
\end{proposition}
\begin{proof}

    Since $\gldim A < \infty$, $\psiforest$ has only one component by \Cref{lemma:infinitegldim}.
   Applying \Cref{lemma:idim:bound:same:tree} implies that
    \[
        \gordim \sAus(A) \leq 6 \height(A) - 2
    \]
    since all indecomposable projective $ \sAus(A) $-modules are isomorphic to $ \pstab{M} $ for some indecomposable non-projective $ A $-module $ M $ (\Cref{theoreminjcoresstableaus}.\eqref{theoreminjcoresstableaus1}). The result follows now from the equality $ \gldim \sAus(A) = \gordim \sAus(A) $, which holds since the former dimension is finite by $ \gldim \sAus(A) \leq 3 \gldim A - 1 $ (see \cite[Proposition II.3.9]{AR}).
   For the optimality of the bound, see \Cref{optimal.bound}.
\end{proof}

\begin{example}\label{optimal.bound}
  Let $ A $ be the Nakayama algebra with coKupisch series $d=[3,3,3,3, \allowbreak 4,5]$, whose coresolution forest is shown in \Cref{fig:coresforest-333345}. 
  Observe that $ \gldim A < \infty $ as $ \psiforest $ is connected. 
  One can verify that $\gordim \sAus(A) = \idim \pstab{S_2} = 10 = 6\height(A) - 2$ by computing $\idim \pstab{M}$ for all indecomposable non-projective $A$-modules $M$ (or only the simple ones by \Cref{theorem:simples}) via the three-distance formula.
\end{example}

  \begin{figure}[t]
    \begin{subfigure}[t]{0.32\textwidth}
      \centering
      \input{figures/fig-coresforest-12344555.tex}
      \caption{$d = [1,2,3,4,4,5,5,5]$}
      \label{fig:coresforest-12344555}
  \end{subfigure}\hfill
    \begin{subfigure}[t]{0.32\textwidth}
      \centering
  \centering
  {
  \newcommand{\nodea}{\node[draw=none, inner sep=2pt] (a) {}
  ;}
  \newcommand{\nodeb}{\node[draw=none, inner sep=2pt] (b) {}
  ;}
  \newcommand{\nodec}{\node (c) {$\overline{2}$}
  ;}
  \newcommand{\noded}{\node (d) {$\overline{1}$}
  ;}
  \newcommand{\nodee}{\node (e) {$0$}
  ;}
  \newcommand{\nodef}{\node (f) {$1$}
  ;}
  \newcommand{\nodeg}{\node (g) {$2$}
  ;}
  \newcommand{\nodeh}{\node (h) {$3$}
  ;}
  \newcommand{\nodei}{\node (i) {$4$}
  ;}
  \newcommand{\nodej}{\node (j) {$5$}
  ;}
  \newcommand{\nodek}{\node (k) {$6$}
  ;}
  \newcommand{\nodel}{\node[draw=none, inner sep=2pt] (l) {}
  ;}
  \newcommand{\nodem}{\node[draw=none, inner sep=2pt] (m) {}
  ;}
  \newcommand{\noden}{\node[draw=none, inner sep=2pt] (n) {}
  ;}
  \begin{tikzpicture}[auto, font=\footnotesize]
  \matrix[column sep=.3cm, row sep=.5cm,ampersand replacement=\&]{
    \nodea  \&  \nodeb  \&           \\
    \nodec  \&  \noded  \&           \\
    \nodee  \&  \nodef  \&  \nodeg  \\
    \nodeh  \&  \nodei  \&           \\
    \nodej  \&  \nodek  \&           \\
    \nodel  \&  \nodem  \&  \noden  \\
  };

  \path[ultra thick]   (d) edge (e)
  	  (d) edge (f)
  	  (d) edge (g)
  	  (e) edge (h)
  	  (f) edge (i)
  	  (h) edge (j)
  	  (i) edge (k);
  \path[ultra thick, dotted]   (a) edge (c)
  	  (b) edge (d)
  	  (k) edge (l)
  	  (k) edge (m)
  	  (k) edge (n);
  \end{tikzpicture}}
      \caption{$d = [1,2,3,3,3,2,2]$}
      \label{fig:coresforest-1233322}
    \end{subfigure}\hfill
    \begin{subfigure}[t]{0.32\textwidth}
      \centering
      \input{figures/fig-coresforest-333345.tex}
      \caption{$d = [3,3,3,3,4,5]$}
      \label{fig:coresforest-333345}
    \end{subfigure}
    \caption{Portions of the coresolution forest $\psiforest$ of the three Nakayama algebras with the displayed coKupisch series.}
    \label{fig:boundsexample2}
  \end{figure}

  \begin{remark}
      Observe that our bound in \Cref{prop:finite:height} is sharper by one than the classical bound in \cite{AR} when $ \findim A $ is even, as
        \[
            6 \height(A) - 2 = 3 \findim A - 2 = 3 \gldim A - 2
        \]
    by \Cref{prop:findimA} in this case.
  \end{remark}

Let us observe via an example that the bound of \Cref{prop:finite:height} can be arbitrarily larger than the actual value of $ \gldim \sAus(A) $, even though it is optimal:

  \begin{example}\label{arbitrarily.larger}
      Let $ A $ be the linear Nakayama algebra with $ n $ vertices and coKupisch series $ d = [1,2,2, \ldots ,2] $. It is straightforward to see that $ \gldim \sAus(A) = 0 $ as all the non-projective $ A $-modules are simple. On the other hand, one computes
      $$ \height(A) = \left\lceil \frac{n-1}{2} \right\rceil ,  $$
      implying that $ 6 \height(A) - 2 \geq 3 n - 8 $.
      We note that
      $ \gldim A = n - 2 $
      by \Cref{prop:findimA}.
  \end{example}

Returning to the case $ \gldim A = \infty $, we prove the main result of this subsection:

\begin{theorem}\label{prop:infinite:height}
    Let $A$ be a Nakayama algebra with $\gldim A = \infty$. Then,
        \[
            \fingordim \sAus(A)  =  \begin{cases}
                6 \height(A) , & \text{if $ \psiforest $ has a non-left vertex} \\
                               & \text{of maximum height;} \\
                6 \height(A)-2 , & \text{otherwise.}
                    \end{cases}
        \]
    Moreover, $ \fingordim \sAus(A) = \idim \pstab{S_i} $, where $ i-1 $ is any rightmost vertex that is non-left and of maximum height, if such a vertex exists, or else it is any rightmost vertex that is left and of maximum height.
\end{theorem}

\begin{proof}
  By \Cref{prop:self_injective,lemma:infinitegldim}, $\psiforest$ has more than one component and at least one leaf.  Assume first that $i-1$ is a rightmost vertex that is right and of maximum height.
  Since $\psiforest$ has more than one component, if $ i $ is in the same component as $ i-1 $, then it must be on the same level as $i-1$ and, therefore, to its right. But then $ i $ is also of maximum height violating the choice of $ i-1 $.

  Now, \Cref{S-lemma:corollaries}.(\ref{lemma:right.vertex}) implies that $\gamma(i)-1$ is
  in the same component of $ \psiforest $ as $i-1$ and one level below.
    Thus, $ \dist(i-1,\gamma(i)-1)=2 \height(A)+1$ by \Cref{cor:distance}, whereas
  $\dist(i-1, i) = \dist(\gamma(i)-1, i) = \infty$.
  Applying
  \Cref{three:distance:formula}, we obtain
  \[
    \idim\pstab{S_i} %
    = 3(\dist(i-1, \gamma(i)-1)-1) %
    = 6\height(A).
  \]

  It remains to show that there is no indecomposable non-projective module $M = \mik[j]$ with $\idim \pstab{S_i} < \idim \pstab{M} < \infty$.
  By \Cref{lemma:idim:bound:same:tree}, we only need to consider modules $M(j,k)$ for which
  exactly two of the three vertices $j-1$,
  $j+k-1$, and $\gamma(j)-1$ lie in the same component.
    If $\idim \pstab{M} > \idim \pstab{S_i}$, then these vertices are $ j-1 $ and $ j+k-1 $ and their distance is $2\height(A) + 1$ by \Cref{three:distance:formula}. In particular, $ j+k-1 $ is on the level below $ j-1 $ by \Cref{cor:distance} and \Cref{S-lemma:corollaries}.(\ref{S-lemma:levels}).
    By the latter corollary, $ \gamma(j)-1 $ is on the same level as $ j+k-1 $ since the level function is weakly increasing. Therefore, $ j+k-1 $ is to the left of $ \gamma(j)-1 $, which in turn is to the left of $ j-1 $ by \Cref{lemma:largest:child:non:leaf}. It follows that $ \gamma(j)-1 $ is in the same component as the other two vertices, a contradiction.

  Assume now that any leaf of height $\height(A)$ is a left
  leaf, and let $i-1$ be a rightmost such leaf.
  Then $\gamma(i)-1$ is in the previous component, because all
  vertices to the left of $i-1$ and on the same floor must be leaves, otherwise the maximality of $ \height(i-1) $ would be violated
  (using \Cref{lemma:largest:child:non:leaf}).
  On the other hand, vertex $i$ is on the same floor as $i-1$ since the spine is to the right of $i-1$.
  Moreover, $i-1\vee i$ is
  on the spine because, otherwise, $\height(i-1)=\height(i)$ would violate the choice of $ i-1 $.
    Thus, $\dist(i-1, i) = 2\height(A) < \infty$ whereas
  $\dist(i, \gamma(i)-1) = \dist(i-1, \gamma(i)-1) = \infty$. Applying
  \Cref{three:distance:formula}, we obtain
  \[
    \idim\pstab{S_i} %
    = 3(\dist(i, i-1)-1)+1 %
    = 6\height(A) - 2.
  \]

  Again, we argue that there is no $M = \mik[j]$ for which  $\idim \pstab{S_i} < \idim \pstab{M} < \infty$.
  As before, we only need to consider the case when exactly two of $j-1$,
  $j + k - 1$, and $\gamma(j)-1$ lie in the same component; let us call this pair $x$ and $y$.
  If $\dist(x,y) \leq 2 \height(A)$, then $\idim \pstab{M} \leq 3(2 \height(A) - 1) + 1 = 6 \height(A) - 2$.
  So $\idim \pstab{M} > \idim \pstab{S_i}$ implies $\dist(x,y) = 2 \height(A) + 1$, which is only possible when $ x $ and $ y $ lie on different levels.
  If we were in Case $1$ or $3$ of \Cref{three:distance:formula}, then $\idim \pstab{M}$ would be odd, contradicting \Cref{prop:even}.
  Therefore, we must be in Case $2$, meaning $x$ and $y$ are $j - 1$ and $\gamma(j) - 1$.
  But then $ j-1 $ has maximum height by \Cref{cor:distance}, and \Cref{lem:non.IG} implies $ \dist(x,y) = \infty $, a contradiction.
\end{proof}

Combining \Cref{prop:infinite:height} with \Cref{prop:findimA}, we get:

\begin{corollary}\label{cor:finGordim&findim}
    Let $ A $ be a Nakayama algebra with $ \gldim A = \infty $. Then
        \[
            \fingordim \sAus(A)  =  \begin{cases}
                3 \findim A , & \text{if $ \, \findim A $ is even;} \\
                3 \findim A + 1 , & \text{otherwise.}
                    \end{cases}
        \]
\end{corollary}
\begin{proof}
    If $ A $ is not self-injective, the equality follows from \Cref{prop:findimA,prop:infinite:height}. If $ A $ is self-injective, then it follows from \Cref{prop:self_injective}.
\end{proof}

When $ A $ is Iwanaga--Gorenstein, there is a straightforward relationship between the Gorenstein dimensions of $ A $ and $ \sAus(A) $:

\begin{proposition}\label{prop:iwanaga-gorenstein}
    For a Nakayama algebra $A$, if $ A $ is Iwanaga--Gorenstein then so is $\sAus (A) $.  If in addition $\gldim A = \infty$, then $$\gordim \sAus (A) = 3 \gordim A.$$
\end{proposition}
\begin{proof}\label{proof:iwanaga-gorenstein}
Since the injective dimensions of all projective $ A $-modules are finite, we have $\idim \pstab M  \leq 3 \idim P(M) < \infty$ for every non-projective $ A $-module $ M $ by \Cref{three:distance:formula}. Therefore, $\gordim \sAus (A) < \infty $ by \Cref{theoreminjcoresstableaus}.(\ref{theoreminjcoresstableaus1}). The same argument about the opposite algebras shows that $ \sAus(A) $ is Iwanaga--Gorenstein.
Moreover, $\gordim \sAus (A) \leq \sup\left\{ 3 \idim P \right\} = 3\gordim A$, where $ P $ ranges over projective $ A $-modules.

We now show the reverse inequality when $ \gldim A = \infty $ by constructing a projective $\sAus (A) $-module with injective dimension equal to $3 \gordim A$.
Suppose $\gordim A$ is achieved at the projective module $P = M(i, c_i)$.
Since $\gordim A < \infty$, $i$ and $i + c_i - 1$ lie in the same component of $ \psiforest $.
Since $\gldim A = \infty$, $\psiforest$ contains more than one component (\Cref{lemma:infinitegldim}), and $ i-1 $ and $ \gamma(i)-1 $ do not lie on the same floor
(\Cref{cor:finite.gldim}).
Therefore, we can pick $k$ so that $i + k - 1$ lies in a different component.
Consider $M = \mik$, which is non-projective by construction, and observe that $P(M) = P$.
By \Cref{three:distance:formula}, we have $\idim{M} = \idim{\Omega(M)} = \infty$ while $\idim{P(M)} < \infty$, so $\idim \pstab{M} = 3\idim{P(M)} = 3\gordim A$.
Therefore $\gordim \sAus(A) \geq 3 \gordim A$, which completes the proof.
\end{proof}

\begin{remark}
    \Cref{prop:iwanaga-gorenstein,cor:finGordim&findim} recover \cite[Corollary 3.26]{Se2}, i.e., that $ \gordim A = \findim A $ is even when $ A $ is an Iwanaga--Gorenstein Nakayama algebra with infinite global dimension.
\end{remark}

\Cref{prop:infinite:height} provides a practical way to compute the finitistic dimension of the stable Auslander algebra $ \sAus(A) $ through the equality with its finitistic Gorenstein dimension when $ \sAus(A) $ is Iwanaga--Gorenstein. However, this is not the case in general as illustrated by the next example:

\begin{example}\label{fGd.less.than.findim}
For the Nakayama algebra $A$ with coKupisch series $d = [4,4,4,3]$,
    \[
        \fingordim \sAus(A) = 4 < 5 = \findim \sAus(A) = \findim \sAus(A)^\op  .
    \]
The last two dimensions were obtained using the GAP-package \cite{QPA} to compute the projective dimensions of all the finitely many indecomposable modules over $ \sAus(A) $ and $ \sAus(A)^\op $.
\end{example}

Next, we apply \Cref{prop:infinite:height} to two concrete algebras:

\begin{example}
Consider the Nakayama algebra $ A $ with coKupisch series $ d = [2, \allowbreak 3, \allowbreak 3,4,4,2] $, whose coresolution forest is shown in \Cref{fig:coresforest-233442}.
This forest has more than one component and so $ \gldim A = \infty $.
It has two vertices (modulo $n=6$) of maximum height $ \height(A)=1 $: vertex $1$, which is right, and vertex $2$, which is left.
It can be seen that $\fingordim \sAus(A) = \idim \pstab{S_2} = 6 = 6 \height(A)$, and $\idim \pstab{S_3} = 4 < \fingordim \sAus(A)$.
\end{example}

\begin{example}
Now consider the Nakayama algebra $ A $ with coKupisch series $ d = [2,3,4,4,5,5,6,5,6,6,3,4,3,4,4,4,4] $, whose coresolution forest is shown in \Cref{fig:coresforest-23445565663434444}.
It has no right vertices of height $\height(A)=2$ and two rightmost left vertices of height $ 2 $, namely vertices $6$ and $14$.
One computes $\fingordim \sAus(A) = \idim \pstab{S_7} = \idim \pstab{S_{15}} = 10 = 6 \height(A) - 2$.
Note that $\fingordim \sAus (A) $ is also attained at $\pstab{M}$ for the non-simple module $M = \mik[6][2]$.
\begin{figure}[t]
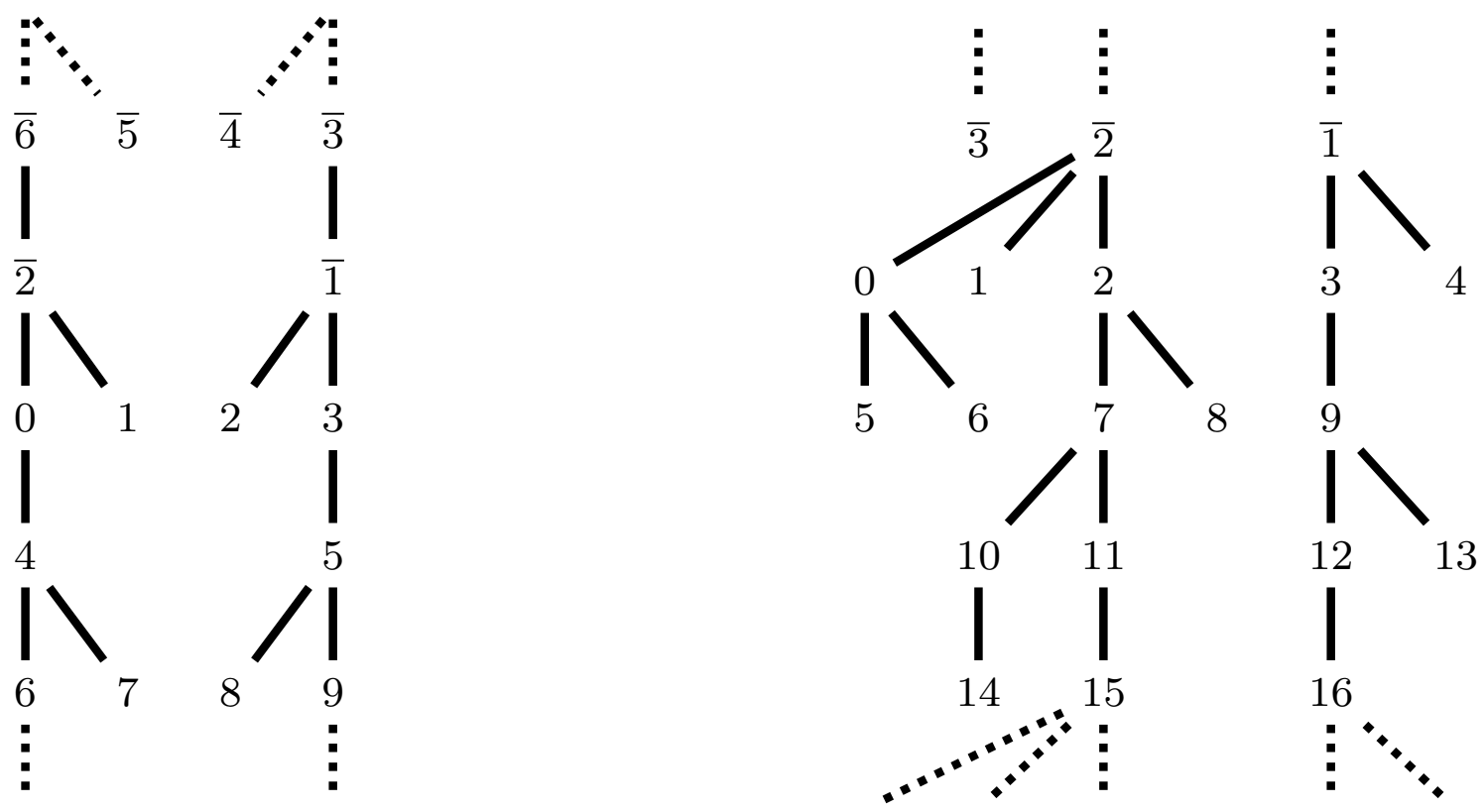

   \begin{subfigure}[t]{0.38\textwidth}
      \centering
        \input{figures/fig-coresforest-233442.tex}
        \caption{$d = [2,3,3,4,4,2]$}
        \label{fig:coresforest-233442}
    \end{subfigure}\hfill
    \begin{subfigure}[t]{0.58\textwidth}
      \centering
      \input{figures/fig-coresforest-23445565663434444.tex}
      \caption{$d = [2,3,4,4,5,5,6,5,6,6,3,4,3,4,4,4,4]$}
      \label{fig:coresforest-23445565663434444}
    \end{subfigure}
    \caption{Portions of the coresolution forests of the two Nakayama algebras with the displayed coKupisch series.}
      \label{fig:boundsexample}
\end{figure}
\end{example}

We conclude this subsection by showing that non-projective simple $ A $-modules give rise to projective $ \sAus(A) $-modules of finite injective dimension, regardless of whether the global dimension of $A$ is finite:

\begin{proposition}\label{thm:idim-simple-finite}
    The injective dimension $ \idim \pstab S $ is finite for every non-pro\-jec\-tive simple module $ S $ over a Nakayama algebra $A$.
\end{proposition}

\begin{proof}
  Consider the simple $A$-module $S_i$.  Then $\idim \pstab{S_i}$
  depends on the distances between all possible pairs of vertices
  among $i-1$, $i$, and $\gamma(i)-1$ in the coresolution forest
  $\psiforest$.  If any of these distances is finite, that is, if any
  two of the vertices are in the same component of $ \psiforest $,
  then $\idim \pstab{S_i}$ is finite.

  If $i - 1$ is not a leaf, $\gamma(i) - 1$ is a child of $i-1$ and
  thus in the same component as $i-1$ by \Cref{lemma:largest:child:non:leaf}.
  If $i-1$ is a leaf, and $i - 1$ and $i$ are in different components,
  then $i-1$ is a right vertex. Therefore, $\gamma(i)-1$ and $i-1$ are in the same component by \Cref{S-lemma:corollaries}.\eqref{lemma:right.vertex}.
\end{proof}

\subsection{Gorenstein dimension and simple modules}\label{sec:main-results-simples}
The goal of this subsection is to prove that $ \fingordim \sAus(A)$ is attained at a simple module even when $ \gldim A < \infty $ (see \Cref{theorem:simples}), completing the proof of \Cref{mainresult2} together with \Cref{prop:infinite:height,thm:idim-simple-finite} from the previous subsection. Since $ \gldim \sAus(A) < \infty $ in this case by \cite[Proposition II.3.9]{AR}, this is essentially a statement about $ \gldim \sAus(A) $.

We begin with two crucial preliminary results for coresolution forests:

\begin{lemma}\label{lemma:minimal-distance}
  Assume that the coresolution forest $\psiforest$ is connected, and
  let $a , b , c \in \psiforest$ be vertices such that $a$ is to the left of $b$, which in turn is to the left of $c$.  Suppose that, either
  \begin{enumerate}[\rm(i)]
  \item\label{order:abc} $\level(a) = \level(b) = \level(c)$ and
    $\mathcal I = \{a,\dots,c\}$, or
  \item\label{order:bca} $\level(a) - 1 = \level(b) = \level(c)$ and
    $\mathcal I = \{b,\dots,a\}$, or
  \item\label{order:cab} $\level(a) = \level(b) = \level(c) + 1$ and
    $\mathcal I = \{c,\dots,b\}$.
  \end{enumerate}
  Let
  \begin{align*}
    \tilde{\mathcal A} &= \{j\in\mathcal I\mid\dist(j, b)\geq\dist(a, b)\text{ and $j$ left of $b$}\},\\
    \tilde{\mathcal C} &= \{j\in\mathcal I\mid\dist(b, j)\geq\dist(b, c)\text{ and $j$ right of $b$}\},\\
    \tilde{\mathcal B} &= \mathcal I\setminus(\tilde{\mathcal A}\cup\tilde{\mathcal C}).
  \end{align*}
  Then $\tilde{\mathcal A}$, $\tilde{\mathcal B}$, and $\tilde{\mathcal C}$ are intervals that partition $\mathcal I$, where $ \tilde{\mathcal A} $ is to the left of $ \tilde{\mathcal B} $ and $ \tilde{\mathcal B} $ is to the left of $ \tilde{\mathcal C} $.
Moreover, $\dist(\tilde x, \tilde y)\geq\dist(x, y)$ for $x, y\in\{a, b, c\}$ and all elements $\tilde a\in\tilde{\mathcal A}$, $\tilde b\in\tilde{\mathcal B}$, and $\tilde c\in\tilde{\mathcal C}$.
\end{lemma}
\begin{proof}
    Let us first show that the vertices in each of the sets $\tilde{\mathcal A}$, $\tilde{\mathcal B}$ and $\tilde{\mathcal C}$ are on the same level and form intervals.  In Case~\ref{order:abc}, all vertices in $\mathcal I$ are on the same level, because $\level$ is weakly increasing by \Cref{thm:naturally}.(\ref{monotone}). \Cref{lemma:common-ancestor} now implies that $\tilde{\mathcal A}$ is an initial interval of $\mathcal I$, and that $\tilde{\mathcal C}$ is a terminal interval of $\mathcal I$.  Thus, also $\tilde{\mathcal B}$ is an interval.

    In the other two cases, the vertices in $\mathcal I$ are on two levels.  We focus on Case~\ref{order:bca}; Case~\ref{order:cab} is very similar.  We show that in Case~\ref{order:bca}, precisely the vertices in $\tilde{\mathcal A}$ are on the same level as $a$.  Indeed, any vertex in $\tilde{\mathcal A}$ is to the left of $b$
    and thus cannot be on the same level as $b$.  Conversely, any vertex $j < a$ on the same level as $a$ is to the left of $a$, and therefore also to the left of $b$. \Cref{lemma:common-ancestor} then states that $\dist(j, b)\geq \dist(a, b)$, so $j\in\tilde{\mathcal A}$.  Since $\tilde{\mathcal A}$ is the subset of vertices on the same level as $a$, it is a terminal interval of $\mathcal I$.  The remaining vertices of $\mathcal I$ are all on the same level, so \Cref{lemma:common-ancestor} implies that $\tilde{\mathcal C}$ and $\tilde{\mathcal B}$ are also intervals.

    It remains to show that $\dist(\tilde x, \tilde y)\geq\dist(x, y)$ for $x,y\in\{a,b,c\}$.  We use that
  \[
  \dist(\tilde x, \tilde y)\geq\dist(x, y) \Longleftrightarrow \level(\tilde x\vee\tilde y)\leq\level(x\vee y)
  \]
    due to \Cref{prop:distance} and the equalities $ \level( \tilde x ) = \level( x )$ and $ \level( \tilde y ) = \level( y )$.

    We begin with the inequality $\level(\tilde a\vee\tilde c)\leq\level(a\vee c)$.  Indeed, applying \Cref{lemma:common-ancestor} twice, and \Cref{prop:distance} to the inequalities $\dist(\tilde a, b)\geq\dist(a,b)$ and $\dist(b, \tilde c)\geq\dist(b,c)$, yields
  \begin{align*}
    \level(\tilde a\vee\tilde c)&=\min\{\level(\tilde a\vee b), \level(b\vee\tilde c)\}\\
    &\leq\min\{\level(a\vee b), \level(b\vee c)\} = \level(a\vee c).
  \end{align*}

    Finally we prove that $\level(\tilde a\vee\tilde b)\leq\level(a\vee b)$, and the proof of $\level(\tilde b\vee\tilde c)\leq\level(b\vee c)$ is very similar.  By definition of $\tilde{\mathcal A}$ and $\tilde{\mathcal B}$ we have
  \begin{align}
    \dist(\tilde a, b)&\geq\dist(a, b),\text{ and}\tag{d1}\label{eq:H1}\\
    \dist(\tilde b, b)&<\dist(a, b)\tag{d2}\label{eq:H2}\text{ or $\tilde b$ is not to the left of $b$.}
  \end{align}
  Applying \Cref{prop:distance} to \eqref{eq:H1} yields
  $\level(\tilde a\vee b)\leq\level(a\vee b)$.  This implies the
  claim if $\tilde b$ and $b$ coincide.

  If $\tilde b$ is to the right of $b$, we apply
  \Cref{lemma:common-ancestor} to $\tilde a$, $b$ and $\tilde b$.  We
  then find that
  $\level(\tilde a\vee\tilde b)\leq\level(\tilde a\vee
  b)\leq\level(a\vee b)$.

  If $\tilde b$ is to the left of $b$, Case~\ref{order:bca} is impossible
  because then all elements of $\mathcal I$ to the left of $b$ are in $ \tilde{ \mathcal A } $.  Thus $a$, $\tilde a$, $b$ and $\tilde b$ are all on the same level.  Applying
  \Cref{prop:distance} to \eqref{eq:H1} and \eqref{eq:H2} yields
  $\level(\tilde a\vee b)\leq\level(a\vee b)<\level(\tilde b\vee b)$.
  Applying \Cref{lemma:common-ancestor} to $\tilde a$, $\tilde b$ and
  $b$ we find
  \[
    \level(\tilde a\vee b) = \min\{\level(\tilde a\vee\tilde b),
    \level(\tilde b\vee b)\},
  \]
  so $ \level(\tilde a\vee\tilde b) = \level(\tilde a\vee b) \leq \level(a\vee b) $.
\end{proof}

    In what follows, we write $ \downarrow \! \mathcal S $ for the downward closure of a subset $\mathcal S \subseteq \psiforest$, i.e., for the set containing $ \mathcal S $ and all vertices that have an ancestor in $ \mathcal S $.

\begin{lemma}\label{lemma:larger-simple}
  Assume the hypothesis of \Cref{lemma:minimal-distance} and consider
  the partition of $ \downarrow \! \mathcal I $ into the three disjoint
  subsets $\mathcal A = \downarrow \! \tilde{\mathcal A} $,
  $\mathcal B = \downarrow \! \tilde{\mathcal B} $ and
  $\mathcal C = \downarrow \! \tilde{\mathcal C} $.
  Assume furthermore that
  $ \downarrow \! \mathcal I $
  contains all vertices greater than $\min\mathcal I$.
    Then, we can choose a vertex~$i'$ such that each of $\mathcal A$,
  $\mathcal B$ and $\mathcal C$ contains precisely one of
 $i'-1$, $i'$
  and $\gamma(i')-1$.
\end{lemma}
\begin{proof}
  Since $\psiforest$ is connected and thus contains exactly one doubly infinite path by \Cref{prop:covering}, precisely two sets among $\mathcal A$, $\mathcal B$ and $\mathcal C$ are finite.

  Let $\mathcal X$, $\mathcal Y$ and $\mathcal Z$ be obtained by
  cyclically permuting $\mathcal A$, $\mathcal B$ and $\mathcal C$,
  such that $\mathcal X$ is the
  leftmost finite set where the level
  of the largest vertex is minimal. Let $ i'-1 $ be the largest vertex in $ \mathcal X $.

  Let $\mathcal S$ be either $\mathcal Y$ or $\mathcal Z$.  We show that
  $\mathcal S\cap[i',\gamma(i')-1]$ is non-empty.  Suppose otherwise, then
  \Cref{S-lemma:corollaries}.(\ref{S-lemma:ancestors})
  implies that $\max\mathcal S$ exists.  If $\mathcal S$ is to
  the right of $\mathcal X$, then $\level(\max\mathcal S) < \ell$ since every vertex on level $ \ell $ is in $ \mathcal X $, to the left of $ \mathcal X $, or in $[i',\gamma(i')-1]$ by \Cref{S-lemma:corollaries}.(\ref{S-lemma:interval}). If
  $\mathcal S$ is to the left of $\mathcal X$, then $\level(\max\mathcal S)<\ell+1$ as every vertex on level $ \ell+1 $ is either to the right of $ \mathcal X $ or in $[i',\gamma(i')-1]$ again by \Cref{S-lemma:corollaries}.(\ref{S-lemma:interval}). In both cases, the assumption that $\mathcal S\cap[i',\gamma(i')-1]$ is empty violates the choice of $ \mathcal X $.

  Since all vertices in $[i',\gamma(i')-1]$ are in $\downarrow\!\mathcal I$, we have $i'\in \mathcal Y\cap[i',\gamma(i')-1]$ and
  $\gamma(i')-1\in\mathcal Z\cap[i',\gamma(i')-1]$.
\end{proof}

We are now ready to prove the main result of this subsection:

\begin{theorem}\label{theorem:simples}
    Let $A$ be a Nakayama algebra with $\gldim A < \infty$.  Then there is a simple non-projective $A$-module $S$ such that $\gldim \sAus(A)= \stabidim[S] < \infty $.
\end{theorem}
\begin{proof}
    Let $M=\mik$ be a non-projective indecomposable $A$-module such that $\stabidim < \infty$.
    The maximum of these dimensions is the finitistic Gorenstein dimension of $\sAus (A)$.

    We construct a simple $A$-module $S=\mik[i'][1]$ such that
    \[
    \stabidim\leq\stabidim[S] < \infty .
    \]

    If $\stabidim=0$ we can choose any simple module.  For
    simplicity, we exclude this case from now on.  Thus, by \Cref{three:distance:formula} and \Cref{lemma:largest:child:non:leaf} we can assume that $i-1$ is a leaf of $\psiforest$.

    Since $i-1$ is a leaf, no two distinct vertices in $\mathcal I=[i-1,\gamma(i)-1]$ are related by ancestorship (\Cref{S-lemma:corollaries}.(\ref{S-lemma:no:ancestors})).
    Let $a$ be the leftmost, $b$ the middle and $c$ the rightmost vertex among $i-1$, $i+k-1$ and $\gamma(i)-1$, and assume the setup of \Cref{lemma:minimal-distance}. We show that $\stabidim\leq\stabidim[S]$ for $S=\mik[i'][1]$, where $i'$ is as in \Cref{lemma:larger-simple}.

    Let $a'\in\mathcal A$,
    $b'\in\mathcal B$ and $c'\in\mathcal C$.  By
    definition, these have unique
    ancestors in $\tilde{\mathcal A}$, $\tilde{\mathcal B}$ and
    $\tilde{\mathcal C}$, respectively, which we denote by
    $ \tilde a $, $ \tilde b $ and $ \tilde c $.  Thus, $ x' \vee y' = \tilde x \vee \tilde y $ for $x, y \in \{a, b, c\}$.  \Cref{lemma:minimal-distance} now implies
    \[
      \dist(x', y') \geq \dist (\tilde x, \tilde y) \geq \dist(x, y).
    \]
    Moreover,
    $\dist(x', y')=\dist(x, y)$ is only possible if both $x'$ and
    $y'$ are in $\mathcal I$.

    As in \Cref{three:distance:formula}, let $d+1$ be the minimal
    distance between any two of the three vertices $i-1$, $i+k-1$ and
    $\gamma(i)-1$.  Let $d'+1$ be the minimal distance between any
    two of the three vertices $i'-1$, $i'$ and $\gamma(i')-1$.  By \Cref{lemma:larger-simple},
    each of the sets $\mathcal A$, $\mathcal B$ and $\mathcal C$
    contains exactly one of $i'-1$, $i'$ and $\gamma(i')-1$.  Thus we
    have $d'\geq d$.  If $d'> d$, it follows that
    $\stabidim[S]\geq 3d'-1\geq 3d+2>\stabidim$.

    If $d'=d$, any
    two vertices among $i'-1$, $i'$ and $\gamma(i')-1$ with minimal distance are in $\mathcal I$.  Since $\min\mathcal I\leq i'-1 < i' < \gamma(i')-1$, both of $ i'-1 $ and $ i' $ must be in $ \mathcal I $.

    If $\gamma(i')-1$ is
    not in $\mathcal I$, then $\idim\Omega(S) = \dist(\gamma(i')-1, i')-1$ and
    $\idim P(S)=\dist(i'-1,\gamma(i')-1)-1$ are larger than $ d = d' $. Therefore,
    $\stabidim[S]=3d+1\geq\stabidim$.

    Otherwise, $i'-1$ and $i-1$ are both in the same set among
    $ \tilde{\mathcal A} $, $ \tilde{\mathcal B} $ and $ \tilde{\mathcal C} $. As the same holds for $ \gamma(i)-1 $ and $ \gamma(i')-1 $, \Cref{lemma:largest:child:non:leaf} implies that
    $\gamma(i')-1 = \gamma(i)-1$.  Because $i'-1$ and $i'$ are in
    different sets, $i'-1$ is the rightmost vertex of the set
    containing $i-1$ and $i'$ is the leftmost vertex of the set
    containing $i+k-1$.

    If $\stabidim=3d-1$, then $\stabidim[S]\geq\stabidim$.  Suppose
    that $\stabidim=3d$.  We have to show that
    $\stabidim[S]\neq 3d-1$, that is, $\dist(i', \gamma(i')-1)$ is
    not among the minimal distances.  This is indeed the case,
    because $\idim P(M)<\idim\Omega(M)$ and
    \Cref{lemma:minimal-distance} implies
    \[
      \dist(\gamma(i)-1, i-1) <
      \dist(i+k-1,\gamma(i)-1)\leq\dist(i',\gamma(i')-1).
    \]

    Finally, suppose that $\stabidim = 3d+1$.  Then
    $\dist(i-1, i+k-1)$ is strictly less than the distances involving
    $\gamma(i)-1$.  By \Cref{lemma:common-ancestor}, it is not
    possible that $a=i+k-1$ and $c=i-1$.  Therefore, $ i-1 \leq i'-1 < i' \leq i+k-1 $ are all on the same level, implying that
    $ d \leq \dist(i'-1, i')\leq\dist(i-1, i+k-1) = d $.
    Combining this with
    \Cref{lemma:minimal-distance}, we obtain
    \begin{align*}
      \dist(i'-1, i')& = \dist(i-1, i+k-1)\\
                     &<\min\{\dist(i+k-1,\gamma(i)-1),\dist(\gamma(i)-1, i-1)\}\\
                     &\leq\min\{\dist(i',\gamma(i')-1),\dist(\gamma(i')-1, i'-1)\},
    \end{align*}
    which yields $ \stabidim[S] = 3d+1 = \stabidim $.
\end{proof}

   We note that in each of the proofs of \Cref{theorem:simples} and \Cref{prop:infinite:height}, the constructed simple $A$-module, say $S = S_i$, is such that $i-1$ is a leaf in the coresolution forest of $A$. In the case where $ \gldim A = \infty $, it holds in addition that $ i-1 $ is a vertex of maximum height.
    
    The next example shows that simple $ A $-modules $ S_i $ such that $ i-1 $ is a leaf of maximum height do not always yield the finitistic Gorenstein dimension of $ \sAus(A) $ when $ \gldim A < \infty $:

\begin{example}
  Let $A$ be the Nakayama algebra with coKupisch series $ d = [1,2,3, \allowbreak 3, 3,2,2] $ (and Kupisch series $ c = [3,3,3,2,2,2,1] $), whose coresolution forest $\psiforest$ is shown in \Cref{fig:coresforest-1233322}.
  Since $ d_0 = 1 $, $ A $ is a linear Nakayama algebra and thus $ \gldim A < \infty $.
    The only vertex in $ \psiforest $ modulo $  n=7 $ where the height of $ A $ is attained is the vertex $ 5 $ ($ \height(A) = \height(5) = 3 $). However, the simple $ A $-module $ S_6 $ is projective ($c_6 = 1$), implying that the simple $ A $-module $ S=S_i $ from \Cref{theorem:simples} is necessarily such that $ i-1 $ is not of maximum height. 
    Observe that $ \gldim \sAus(A) = \idim \pstab{S_3} $ by the discussion preceding this example, since $ 2 $ is the only leaf modulo $ 7 $ such that the simple $ A $-module at $ 2+1 $ is not projective.
\end{example}

\subsection{Gorenstein symmetry conjecture}\label{sec:main-results-symmetry}
This subsection is devoted to the proof of \Cref{mainresult3}, i.e., the validity of the Gorenstein Symmetry Conjecture for stable Auslander algebras of Nakayama algebras; see \Cref{proposition:Gorenstein:symmetry}. We also prove the analogous left-right symmetry for the finitistic Gorenstein dimension of these algebras; see \Cref{fingordim.symmetry}.

Since the Gorenstein dimension of an algebra and its opposite are equal whenever they are both finite (\cite[Lemma A]{Za}), \Cref{mainresult3} is equivalent to:
\begin{proposition}
    \label{proposition:Gorenstein:symmetry}
    For every Nakayama algebra $A$,
        \[
            \gordim \sAus (A) = \infty \, \iff \, \gordim \sAus (A)^\op = \infty .
        \]
\end{proposition}

We begin by interpreting $\gordim \sAus (A) = \infty$ combinatorially in the coresolution forest of $ A $.

\begin{theorem}
    \label{theo:three:components:left:vertex}
    For every Nakayama algebra $A$, $\gordim \sAus(A) = \infty$ if and only if its coresolution forest $\psiforest$ has at least three components and a left vertex.
\end{theorem}
\begin{proof}
      By the three-distance formula (\Cref{three:distance:formula}), $\gordim \sAus (A) = \infty$ if and only if there exist $i$ and $1 \leq k < c_i$ such that $i-1, i + k -1$ and $ \gamma(i)-1 $ are all in different components of $\psiforest$.
      Recall that $\kupisch$ denotes the Kupisch series of $ A $.

    If $\psiforest$ has at least three components and a left vertex, then we can choose $ i-1 $ to be a leftmost vertex that is left and of maximum height among left vertices. Note that $i-1 $ is necessarily a leaf and a leftmost vertex on its floor. Since $ \psiforest $ is not connected, \mbox{$\gamma(i)-1$} lies in the previous component of $ i-1 $ and a level below $ i-1 $  by \Cref{lemma:largest:child:non:leaf}.
    Therefore, $k$ can be chosen so that $i + k - 1$ lies in any component other than the ones containing $ i-1 $ or $ \gamma(i)-1 $.

    Conversely, if $ \gordim \sAus(A) = \infty $ then $ \psiforest $ has at least three components due to the three-distance formula. Furthermore, if every vertex of $ \psiforest $ is right or a vertebra, then $ i-1 $ and $\gamma(i)-1 $ are in the same component of $ \psiforest $ for every $ i $ due to \Cref{S-lemma:corollaries}.(\ref{lemma:right.vertex}), a contradiction.
\end{proof}

\begin{figure}[t]
    \centering
    \input{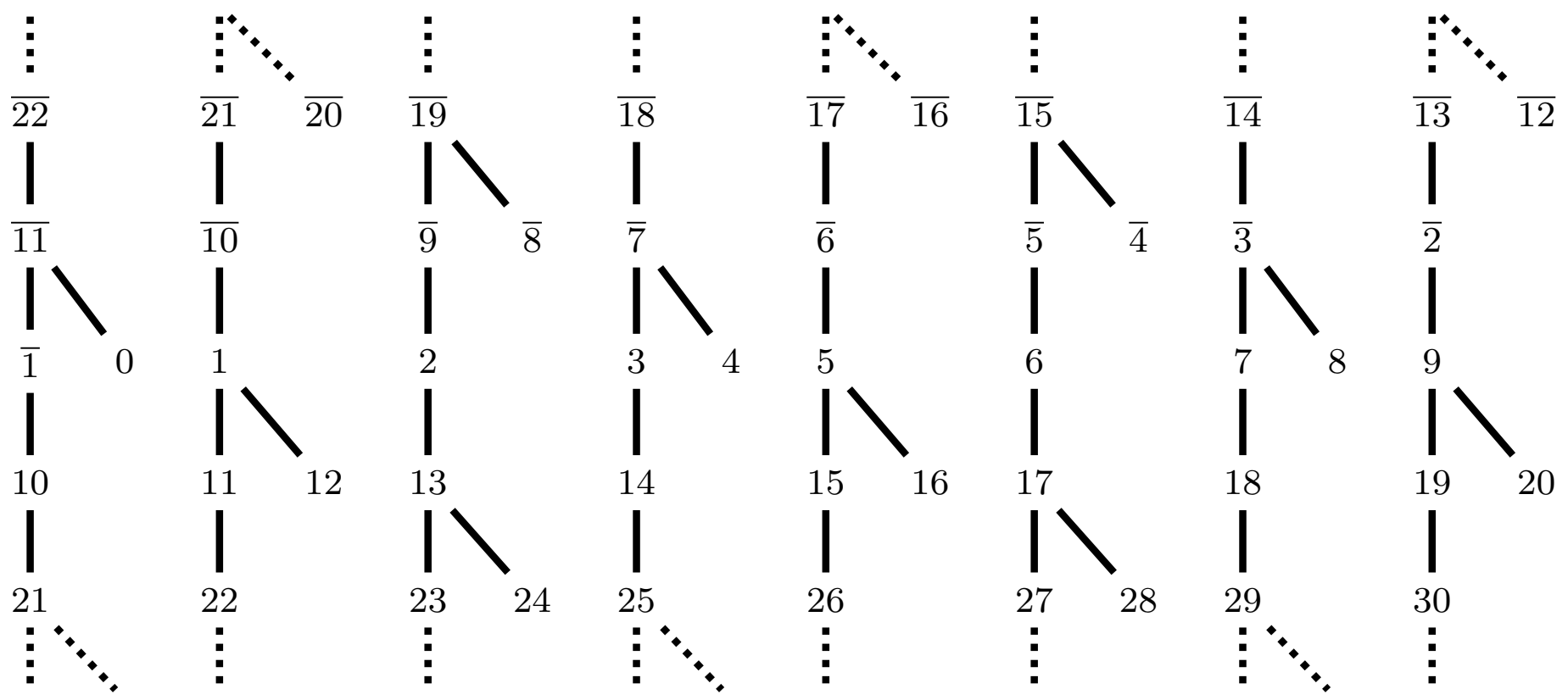}
    \caption{The coresolution forest $\psiforest$ of the Nakayama algebra with coKupisch series $ d=[11,11,11,10] $; it has $8$ components and no left vertex.}
    \label{fig:coresforest-11111110}
    \end{figure}

The following example illustrates that the second requirement in \Cref{theo:three:components:left:vertex} for $ \gordim \sAus(A) $ to be infinite is not redundant:
\begin{example}
    Let $A$ be the Nakayama algebra with coKupisch series $ d = [11,11, \allowbreak 11,10] $.
    Then, its coresolution forest $\psiforest$ has $8$ components as shown in \Cref{fig:coresforest-11111110}, but no left vertex. One may compute that $\gordim \sAus(A) = 6 < \infty $ by \Cref{theo:three:components:left:vertex,prop:infinite:height}, attained at $\pstab{S_1}$.
\end{example}

The opposite algebra $A^\op$ of a Nakayama algebra $A$ is also a Nakayama algebra.  Therefore, to prove \Cref{mainresult3}, we need to connect the coresolution forests 
of $A$ and $A^\op$ in view of \Cref{theo:three:components:left:vertex}.

For any bound quiver algebra $ A = \K Q / I $, $ A^\op $ can be realized as $ \K Q^\op / I^\op $, where $ Q^\op $ is the quiver obtained from $ Q $ by reversing all arrows and $ I^\op $ is the ideal of $ \K Q^\op $ obtained from $ I $ by reversing all arrows in its elements.

When $ A = \K Q / I $ is a Nakayama algebra, the above presentation of $ A^\op $ does not satisfy the convention that the source of an arrow is smaller than its target, except for the arrow $ n-1 \to 0 $ in the cyclic case. Therefore, we represent $ A^\op $ as $ \K Q / \tilde{I} $, where $ \tilde{I} $ is the ideal that results from $ I $ by replacing every arrow $ a_i \colon i-1 \to i $, where vertex labels are reduced modulo $ n $ in the cyclic case, with $ a_{n+1-i}$.

Under this convention, we have $ c^\op_i = d_{n-i} $ and $ d^\op_i = c_{n-i} $, where $ c^\op_i $ and $ d^\op_i $ denote the values in the Kupisch and coKupisch series of $ A^\op $, respectively. Denoting the coresolution and resolution functions of $ A^\op $ by $ \psi^\op $ and $ \gamma^\op $, respectively, it follows that $ \psi^\op(i)= -\gamma(-i) $ and $ \gamma^\op(i) = -\psi(-i) $ for all $ i \in \Z $.
Therefore:
\begin{lemma}\label{lemma:res:forest}
    For every Nakayama algebra $A$, $ \psiopforest = \overline{\gammaforest} $, where $ \overline{\gammaforest} $ is the digraph that results from $ \gammaforest $ by replacing every arrow $ i \to j $ with an arrow $ -i \to -j $.
\end{lemma}

\begin{remark}\label{minus.forest}
\Cref{lemma:res:forest} establishes a clear visual relationship between the coresolution forest $\psiopforest$ of $A^\op$ and the resolution forest $\gammaforest$ of a Nakayama algebra $A$. Specifically, starting with $\gammaforest$ equipped with a natural labeling as in \Cref{levels.of.res.forests}, one obtains $\psiopforest$ with a natural labeling as in \Cref{thm:naturally} simply by replacing every vertex label $i$ with $-i$. See also \Cref{example:four-forests}.
\end{remark}

We now turn to the relationship between the coresolution forests of a Nakayama algebra and its opposite. Crucial to this is the following observation, which expresses $ \gamma $ in terms of $ \psiforest $, and is an immediate consequence of \Cref{lemma:largest:child:non:leaf}.

 \begin{corollary}\label{S-cor}
     For every $ i \in \Z $, $ \gamma(i) $ is the smallest child of the smallest non-leaf larger than or equal to $i$ in $\psiforest$.
 \end{corollary}

\begin{lemma}
    \label{lemma:cyclic:vertices}
    For any Nakayama algebra $A$ and $ v \in \Z $, we have:
    \begin{enumerate}
        \item\label{1} $ v $ is a $ \psi^\op $-vertebra if and only if $ - v $ is a leftmost $ \psi $-vertex.
        \item\label{2} $ v $ is a leftmost $ \psi^\op $-vertex if and only if $ -v $ is a $ \psi $-vertebra.
        \item\label{3} $ v $ is a non-leftmost left (resp., right) $ \psi^\op $-vertex if and only if $ -v $ is a non-leftmost left (resp., right) $ \psi $-vertex.
    \end{enumerate}
\end{lemma}
\begin{proof}
    For any component $ T $ in $ \psiforest $, let $ i^{(j)} $ denote the leftmost vertex of $ T $ on level $ j \in \Z $. Then, $ \psi(i^{(j)}) $ is the smallest child of the smallest non-leaf greater than or equal to $ i^{(j-1)} $, implying that $ \gamma(i^{(j-1)}) = i^{(j)} $ by \Cref{S-cor}. We deduce that the leftmost vertices of $ T $ form a spine in $ \gammaforest $. By \Cref{lemma:res:forest}, the vertices $ -i^{(j)} $ for $ j \in \Z $ form a spine in $ \psiopforest $. In particular, $ \psiopforest $
    has at least as many components as $ \psiforest $ due to \Cref{prop:covering}, and it remains to show that they have exactly the same number of components for the first equality.

    By the symmetric argument, the negatives of the leftmost vertices in $ \psiopforest $
    constitute a spine in $\psiforest$, and thus the two forests have the same number of components, completing the proof of items \eqref{1} and \eqref{2}.

    For item \eqref{3}, if $ v_1 < v_2 $ are two consecutive $ \psi $-vertebrae, meaning that no $ i $ with $ v_1 < i < v_2 $ is a $ \psi $-vertebra, then $ -v_2 < -v_1 $ are consecutive leftmost $ \psi^\op $-vertices by item \eqref{2}.
    Further, by item \eqref{1}, if $ l $ is the unique leftmost $ \psi $-vertex with $ v_1 < l \leq v_2 $, then $ -l $ is the unique $ \psi^\op $-vertebra with $ -v_2 \leq -l < -v_1 $. Therefore, in the $ \psi^\op $-floor $ [ -v_2 , -v_1-1 ] $, the non-leftmost left vertices form the interval $ [ -v_2+1 , -l-1 ] $ and the right vertices form the interval $ [ -l+1 , -v_1-1 ] $. The former correspond to the non-leftmost left vertices of the $ \psi $-floor containing $ l $ and $ v_2 $, while the latter correspond to the right vertices of the $ \psi $-floor containing $ v_1 $, completing the proof.
\end{proof}

\begin{corollary}
    \label{lemma:components}
    For any Nakayama algebra $A$, all four forests $\psiforest$, $\gammaforest$, $\psiopforest$ and $\gammaopforest$ have the same number of components.
\end{corollary}
\begin{proof}
    The forests $ \psiforest $ and $ \gammaforest $ have the same number of components, as do $ \psiopforest $ and $ \gammaopforest $, by \Cref{lemma:res:forest}.
    The same holds for the forests $ \psiforest $ and $ \psiopforest $ by the proof of \Cref{lemma:cyclic:vertices}.
\end{proof}

This finally allows us to prove the main result of this subsection:

\begin{proof}[Proof of \Cref{proposition:Gorenstein:symmetry}]
    By \Cref{theo:three:components:left:vertex}, it suffices to show that $\psiforest$ has at least three components and a left vertex if and only if $\psiopforest$ does.
    The first part is true by \Cref{lemma:components}.
    Suppose that $\psiforest$ has a left vertex, and let $ l $ be a leftmost such vertex.
    Let $ v > l $ be the unique vertebra on the floor of $ l $.
    Then, by \Cref{lemma:cyclic:vertices} and its proof, $-l$ is a $ \psi^\op $-vertebra and $-v$ is the leftmost vertex in the same $ \psi^\op $-floor as $-l$.
    It follows that $ -v $ is a left $ \psi^\op $-vertex. The converse implication follows by applying the same arguments to $A^\op$ in place of $A$.
\end{proof}

We illustrate the results in this section on an example.

\begin{example}\label{example:four-forests}
We look at the four forests in \Cref{fig:fourforests-4545656456} corresponding to the Nakayama algebra $ A $ with coKupisch series $ d = [4,5,4,5,6,5,6,4,5,6] $.
\begin{figure}[t]
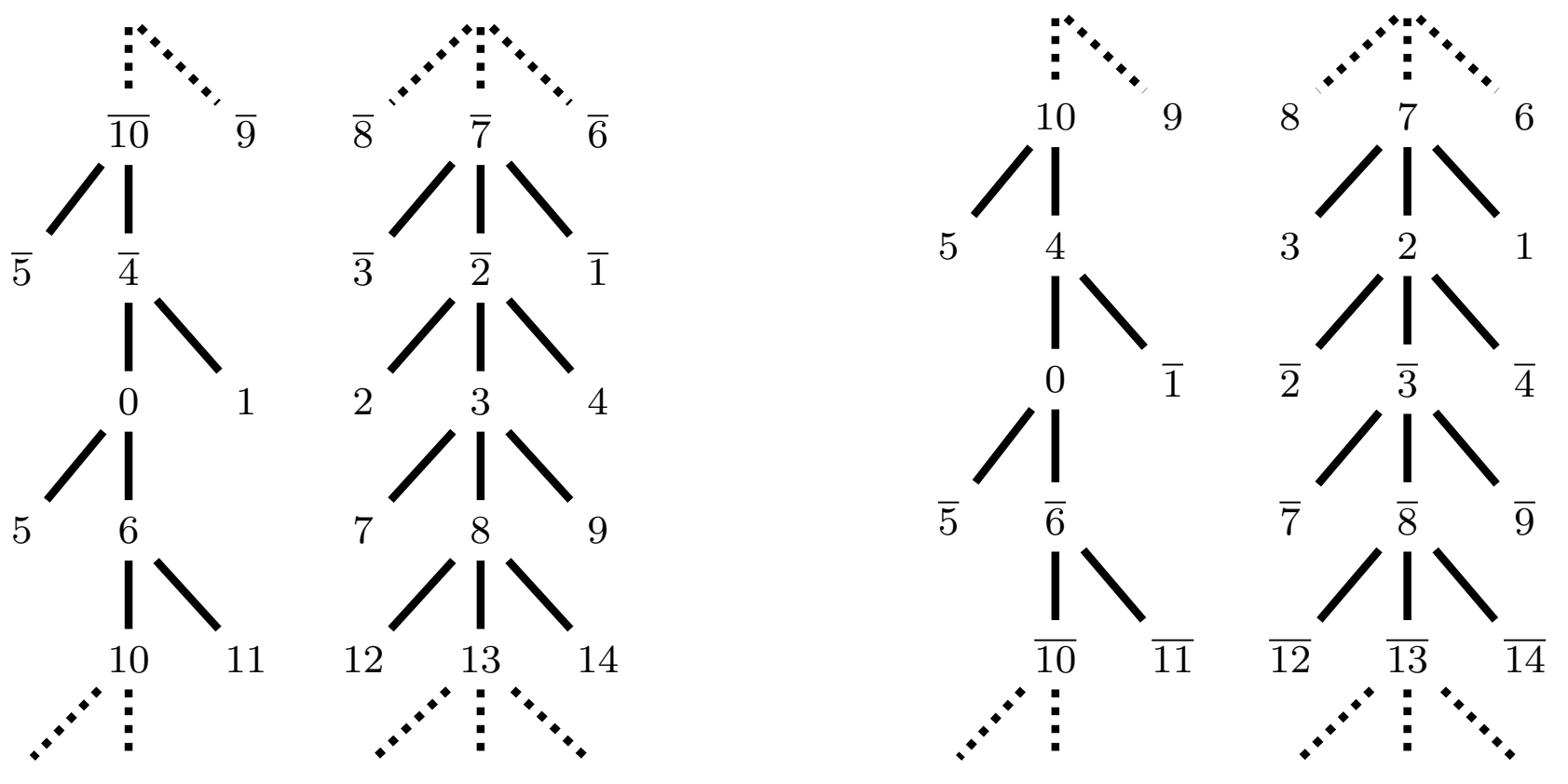
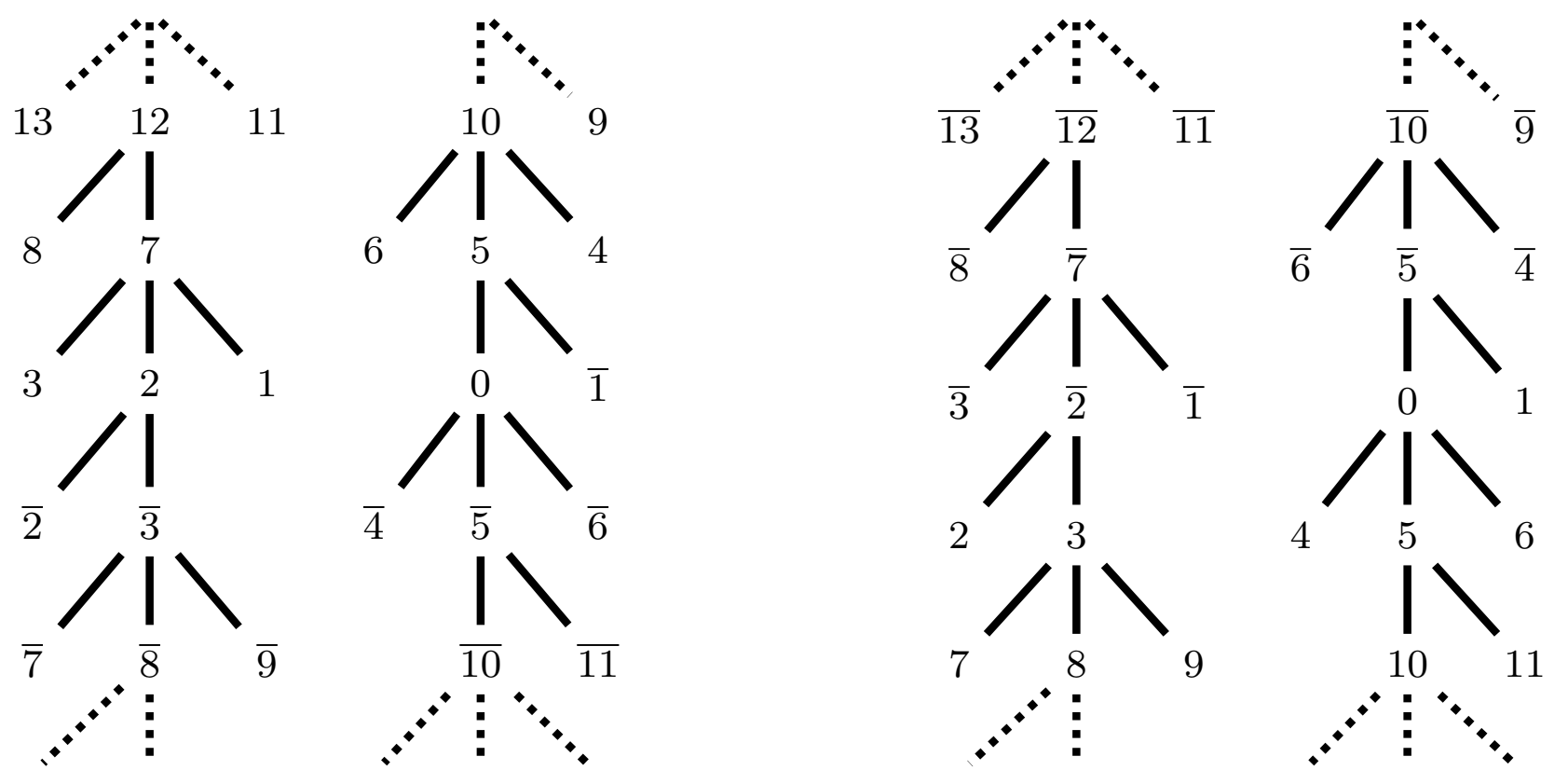

  \centering
  \begin{subfigure}[t]{0.48\textwidth}
    \centering
    \input{figures/fig-coresforest-4545656456}
    \caption{$\psiforest$}
    \label{fig:coresforest-4545656456}
  \end{subfigure}
  \hfill
  \begin{subfigure}[t]{0.48\textwidth}
    \centering
    \input{figures/fig-gammaopforest-4545656456}
    \caption{$\gammaopforest$}
  \end{subfigure}

  \medskip

  \begin{subfigure}[t]{0.48\textwidth}
    \centering
    \input{figures/fig-gammaforest-4545656456}
    \caption{$\gammaforest$}
  \end{subfigure}
  \hfill
  \begin{subfigure}[t]{0.48\textwidth}
    \centering
    \input{figures/fig-psiopforest-4545656456}
    \caption{$\psiopforest$}
  \end{subfigure}
  \caption{The four forests for the Nakayama algebra $ A $ with coKupisch series $ d = [4,5,4,5,6,5,6,4,5,6] $.}
  \label{fig:fourforests-4545656456}
\end{figure}
One observes that the left-right forest pairs are the same with signs inverted, and that the vertebrae in each forest are the leftmost vertices in the corresponding dual forest.
Further, all four forests have the same number of components.
\end{example}

We close this section by proving the analogous left-right symmetry for the finitistic Gorenstein dimension. This is a nontrivial result, since $ \gordim A = \infty = \gordim A^\op $ does not guarantee the equality of the respective finitistic Gorenstein dimensions; see \Cref{fingordim.nonsymmetry}.

\begin{theorem}\label{fingordim.symmetry}
    For every Nakayama algebra $ A $,
        \[
            \fingordim \sAus(A) = \fingordim \sAus(A)^\op   .
        \]
\end{theorem}

\begin{proof}
    In view of \Cref{prop:infinite:height} and the equality $ \height (A) = \height (A^\op) $ (see \cite[Corollary~2, Appendix A.1]{Ri2}), it suffices to show that $ \psiopforest $ has a vertex of maximum height which is a right vertex whenever the same holds for $ \psiforest $, since the converse implication follows by substituting $A^\op$ for $A$.

    Let $ v $ be a $ \psi $-vertex of maximum height $ h = \height(A) $, and suppose it is a right vertex. Then, there is a directed path $ p : v = v_0 \to v_1 \to \ldots \to v_{h-1} \to v_h $ in $ \psiforest $, where all vertices are right except for $ v_h $ which is a $ \psi $-vertebra.
    By \Cref{lemma:cyclic:vertices}.\eqref{3}, $ -v_{h-1} $ is a right $ \psi^\op $-vertex.
    We show that $ -v_{h-1} $ has maximum height $ h $ in $ \psiopforest $.

    By \Cref{S-cor}, $ \gamma(v_{h-1}) $ is the smallest $ \psi $-child of $ v_{h-1} $, and thus $ \gamma(v_{h-1}) \leq v_{h-2} $. Observe that $ \gamma(v_{h-1}) $ cannot be a $ \psi $-vertebra as $ \psi ( \gamma(v_{h-1}) ) = v_{h-1} $ and $ v_{h-1} $ is not a $ \psi $-vertebra. 
    Similarly, $ \gamma^2(v_{h-1}) $ is the smallest $ \psi $-child of the smallest non-$ \psi $-leaf $ w $ that is greater than or equal to $ \gamma(v_{h-1}) $. Since $ v_{h-2} $ is a non-$\psi$-leaf, we have $ \gamma(v_{h-1}) \leq w \leq v_{h-2} $. Therefore, $ \gamma^2(v_{h-1}) $ lies between $ v_{h-3} $ and the vertebra on the same $ \psi $-floor and is not a $ \psi $-vertebra.
    Inductively, $ \gamma^{h-1}(v_{h-1}) $ lies between $ v = v_0 $ and the vertebra on the same $ \psi $-floor and is not a $ \psi $-vertebra; in particular, $ \gamma^{h-1}(v_{h-1}) $ is a right $ \psi $-vertex.

    By \Cref{lemma:res:forest}, we have a directed path
        \[
            -v_{h-1} \to - \gamma(v_{h-1}) \to \ldots \to - \gamma^{(h-1)}(v_{h-1})
        \]
    in $ \psiopforest $ and, by \Cref{lemma:cyclic:vertices}.\eqref{1}, it contains no $ \psi^\op $-vertebra. Therefore, the $ \psi^\op $-height of $ -v_{h-1} $ is at least $ h $. Since it is also at most $ h = \height(A^\op) $, the proof is complete.
\end{proof}

As a consequence, we also recover the left-right symmetry of the finitistic dimension for Nakayama algebras:

\begin{corollary}[\mbox{\cite[Theorem~A.(vii)]{Se2}}]\label{cor:findim.findimop}
    For every Nakayama algebra $ A $,
        \[
            \findim A = \findim A^\op   .
        \]
\end{corollary}

\begin{proof}
    By the proof of \Cref{fingordim.symmetry}, the maximum height $ \height(A) $ in $ \psiforest $ is attained by a right vertex if and only if the maximum height $ \height(A^\op) $ in $ \psiopforest $ is attained by a right vertex. Since $ \height(A) = \height(A^\op) $ by \cite[Corollary~2, Appendix A.1]{Ri2}, the result follows from \Cref{prop:findimA}.
\end{proof}

\section{Further questions and examples}\label{sec:further-questions}

In this section, we ask which of the results we have established so far hold for general representation-finite algebras and provide examples demonstrating the limitations to such potential generalizations.
For the rest of the section, let $A$ be a general non-semisimple representation-finite algebra. As always, we write $\sAus(A) = \underline{\End}_A(X)$ to denote the stable Auslander algebra of $ A $, where $X$ is an additive generator of $\mod A$.

We begin by observing that the finitistic Gorenstein dimension of $ \sAus(A) $ is always finite when $ A $ is not self-injective, as the supremum defining $ \fingordim \sAus(A) $ is in fact a maximum. Indeed, it suffices to consider the injective dimensions of the finitely many indecomposable projective $ \sAus(A) $-modules since
they generate all projective $ \sAus(A) $-modules via direct sums. Furthermore, if $ I $ is an indecomposable non-projective injective $ A $-module, then $ \idim \pstab{I}\leq 1 $ because $\Ext_A^1(X,I)=0$ (see \Cref{theoreminjcoresstableaus}.\eqref{theoreminjcoresstableaus3}), ensuring that the set inside the supremum is nonempty.

The first open question, motivated by \Cref{mainresult2}, is:

\begin{question}
    Can $ \fingordim \sAus(A) $ be attained by computing $ \idim \pstab M$ for a nice subclass of indecomposable non-projective $A$-modules $ M $, for every algebra $ A $ of finite re\-pre\-sen\-ta\-tion type?
\end{question}

The next example shows that, in general, the projective $ \sAus(A) $-modules corresponding to simple $ A $-modules are not enough, unlike the Nakayama case:

\begin{example}\label{not.at.simples}
Let $Q$ be the quiver
\[
\begin{tikzcd}[column sep=large]
1 \arrow[r] & 2 & 3 \arrow[l] \arrow[r] & 4
\end{tikzcd}
\]
and let $A=\K Q$. Since $Q$ is a Dynkin quiver of type $\DynA_4$, the algebra
$A$ is representation-finite and we will use some standard facts about representation-finite hereditary algebras, see for example \cite[Chapter VII]{ASS}. Also note that $A \cong A^{\op}$.

For $1\leq i\leq j\leq 4$, let $M_{ij}$ denote the indecomposable $A$-module
supported on the vertices $i,i+1,\ldots,j$, with all maps between nonzero
vector spaces equal to the identity. The indecomposable projective
$A$-modules are
\[
e_1A=M_{12},\qquad e_2A=M_{22}=S_2,\qquad
e_3A=M_{24},\qquad e_4A=M_{44}=S_4.
\]
Thus the indecomposable non-projective $A$-modules are
\[
\begin{aligned}
&S_1=M_{11},\qquad U=M_{13},\qquad V=M_{14},\\
&W=M_{23},\qquad S_3=M_{33},\qquad Z=M_{34}.
\end{aligned}
\]
Since the projective summands vanish in the stable category, we may take $ \sAus(A)=\underline{\End}_A(X') $, where
\[
X'=S_1\oplus U\oplus V\oplus W\oplus S_3\oplus Z.
\]
Since for any hereditary algebra no nonzero map between two summands of $X'$ factors through a projective $A$-module, we have $ \sAus(A) =\End_A(X')$.
Tedious computations then show that
\[
\sAus(A)\cong \K Q^* /\langle\alpha\varepsilon,\,
\beta\delta-\gamma\zeta\rangle,
\]
where $ Q^* $ is the quiver
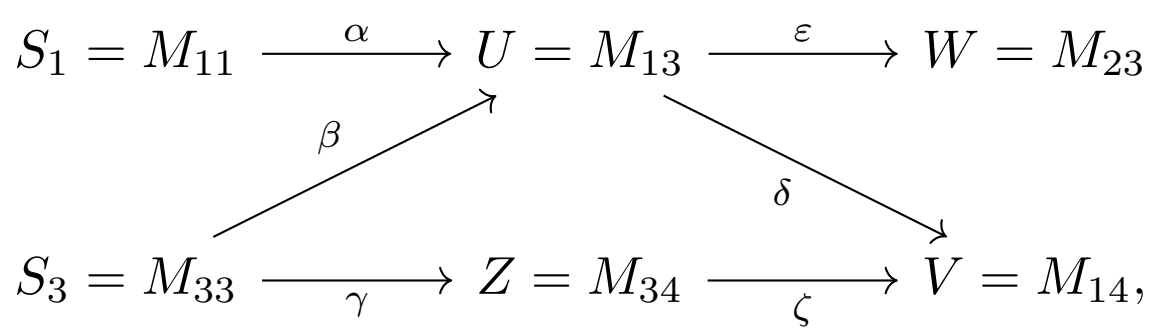
\[
\begin{tikzcd}[column sep=large,row sep=large]
S_1=M_{11} \arrow[r,"\alpha"] & U=M_{13}
  \arrow[r,"\varepsilon"] \arrow[dr,"\delta"'] & W=M_{23} \\
S_3=M_{33} \arrow[ur,"\beta"] \arrow[r,"\gamma"'] & Z=M_{34}
  \arrow[r,"\zeta"'] & V=M_{14} ,
\end{tikzcd}
\]
with vertices labeled by the corresponding
indecomposable $A$-modules.

The injective dimensions of the indecomposable projective $ \sAus(A) $-modules $\pstab M$ are given in the following table:
\[
\begin{array}{c|cccccc}
M&S_1&U&V&W&S_3&Z\\ \hline
\idim\pstab M&1&1&2&2&1&1
\end{array} .
\]
As these are all the indecomposable projective right $\sAus (A) $-modules (see \Cref{theoreminjcoresstableaus}.\eqref{theoreminjcoresstableaus1}), we have
\[
    \gordim \sAus (A) = \fingordim \sAus(A) =2.
\]
The only nonzero modules $\pstab S$ with $S$ simple are $\pstab{S_1}$ and
$\pstab{S_3}$, since $S_2$ and $S_4$ are projective $A$-modules. Hence,
\[
\sup \bigl\{\idim\pstab S\mid S \in \mod A  \text{ is simple} \bigr\}
=1<2=\fingordim\sAus(A).
\]
\end{example}

A second natural question to ask, motivated by \Cref{fingordim.symmetry}, is the following:

\begin{question}\label{Q2}
For which algebras $ A $ of finite representation type beyond Nakayama algebras does $ \fingordim \sAus(A) = \fingordim \sAus(A)^\op $ hold?
\end{question}

\Cref{Q2} is nontrivial since the finitistic Gorenstein dimension of an arbitrary finite-dimensional algebra and its opposite do not coincide in general. We illustrate this with an example:
\begin{example}\label{fingordim.nonsymmetry}
Let $A=\K Q/I$, where $Q$ is the quiver
\[\begin{tikzcd}
	1 & 2 & 3
	\arrow["\alpha", curve={height=-12pt}, from=1-1, to=1-2]
	\arrow["\beta", curve={height=-12pt}, from=1-2, to=1-1]
	\arrow["\gamma", from=1-2, to=1-3]
\end{tikzcd}\]
and $I=\langle \alpha \beta , \beta \alpha , \alpha \gamma \rangle$.
Then, the injective dimensions of the indecomposable projective $A$-modules are given by $0,\infty,\infty$, while the injective dimensions of the indecomposable projective $A^{\op}$-modules are given by $1,0$ and $\infty$. Thus,
$$\fingordim A=0 \neq 1 = \fingordim A^{\op} .$$ We leave the elementary verifications to the reader.
\end{example}

\Cref{mainresult3} motivates the study of the following particular case of the Gorenstein Symmetry Conjecture:

\begin{question}
    Does $ \gordim \sAus(A) = \gordim \sAus(A)^\op $ hold for every algebra $ A $ of finite representation type?
\end{question}

Finally, \Cref{mainresult1} and \Cref{lettercor} prompt the following question regarding the bounds of the Gorenstein dimension:

\begin{question}\label{Q6}
    Can $ \fingordim \sAus(A) $ be bounded in terms of homological data of $A$, such as $ \findim A $ or $ \gordim A $, for every algebra $ A $ of finite representation type?
\end{question}

We note that \Cref{Q6} has a positive answer when $ \gldim A < \infty $, since $ \gldim \sAus(A) \leq 3 \gldim A - 1 $ by \cite[Proposition II.3.9]{AR}.

\section*{Acknowledgements}
This research started at the Junior Research Retreat in Bad Boll 2025 which was supported by the Hausdorff School for Mathematics (HSM) and funded by the Deutsche Forschungsgemeinschaft (DFG, German Research Foundation) under Germany's Excellence Strategy -- EXC-2047/2 -- 390685813.

Odysseas Giatagantzidis is grateful to Steffen Koenig for his generous hospitality during a one-year postdoctoral visit at the University of Stuttgart.
Nupur Jain was funded by the DFG priority program \emph{Combinatorial Synergies} (SPP 2458) --- 539866596.

This project benefited from the use of
the GAP-package \cite{QPA}.

The figures of the forests in this paper were generated using code written by Claude Opus 5.5. 
The other figures were also generated by or formatted with the help of this model.
The authors used Google Gemini Pro for polishing the language and ChatGPT (GPT 6) and Claude Opus 5.5 for proofreading. 
\Cref{not.at.simples,fingordim.nonsymmetry} were generated by ChatGPT (GPT 6).

\end{document}